\documentclass[12pt,psfig]{article}
\usepackage{graphicx,epsfig}
\usepackage{mathrsfs}
\usepackage{amssymb}
\def\al{\alpha}

\def\ld{\lambda}

\def\si{\sigma}

\def\le{\leq}
\def\ge{\geq}

\def\bbb{\begin{eqnarray*}}

\def\eee{\end{eqnarray*}}

\begin{document}

\baselineskip=17pt

\begin{center}

\vspace{-0.6 in} {\large \bf  Regular approximations of spectra of singular discrete\\
linear Hamiltonian systems with one singular endpoint$^*$}
\\ [0.3in]

Yan Liu, Yuming Shi $^{**}$\\

School of Mathematics, Shandong University\\

Jinan, Shandong 250100, P. R. China

\footnote{$^{*}$ This research was partially supported by the NNSF
of China (Grant 11571202) and the China Scholarship Council (Grant 201406220019).}
\footnote{$^{**}$ The corresponding author.}
\footnote{Email addresses: yanliumaths@126.com (Y. Liu), ymshi@sdu.edu.cn (Y. Shi)}
\end{center}

{\bf Abstract.} This paper is concerned with regular approximations of spectra
of singular discrete linear Hamiltonian systems with one singular endpoint.
For any given self-adjoint subspace extension (SSE) of the corresponding
minimal subspace, its spectrum can be approximated by eigenvalues of a sequence of induced regular SSEs,
generated by the same difference expression on smaller finite intervals.
It is shown that every SSE of the minimal subspace has a pure discrete spectrum,
and the $k$-th eigenvalue of any given SSE
is exactly the limit of the $k$-th eigenvalues of the induced regular SSEs; that is, spectral exactness holds,
in the limit circle case. Furthermore, error estimates for the
approximations of eigenvalues are given in this case.
In addition, in the limit point and intermediate cases, spectral inclusive holds.
\medskip

\noindent{\bf 2010 AMS Classification}: 39A10, 47A06, 47A10, 41A99, 34B27.
\medskip

\medskip
\parindent=10pt

\section{ Introduction }

Consider the following discrete linear Hamiltonian system:\vspace{-0.2cm}
$$J\Delta y(t)=(P(t)+\lambda W(t))R(y)(t),\;\;t\in \mathcal{I},      \vspace{-0.2cm}                                 \eqno (1.1_\lambda)$$
where $\mathcal{I}:=\{t\}_{t=a}^{+\infty}$ is an integer interval, $a$ is an integer;
$J$ is the $2n\times 2n$ canonical symplectic matrix, i.e.,
$$
J=\left(\begin{array}{cc} 0&-I_n\\I_n&0\end{array}\right),
$$
with the $n\times n$ identity matrix $I_n$;
$\Delta$ is the forward difference operator, i.e., $\Delta y(t)=y(t+1)-y(t)$;
the weight function $W(t)={\rm diag }\{W_1(t),W_2(t)\}$,
$W_1(t)$ and $W_2(t)$ are $n\times n$ positive semi-definite matrices;
$P(t)$ is a $2n\times 2n$ Hermitian matrix; the partial right shift
operator $R(y)(t)=(y_1^T(t+1), y_2^T(t))^T$ with $y(t)=(y_1^T(t),
y_2^T(t))^T$ and $y_1(t)$, $y_2(t)\in \mathbf{C}^n$;  $\lambda$ is a complex spectral parameter.

It is evident that $P(t)$ can be blocked as
$$
P(t)=\left(\begin{array}{cc}
       -C(t)& A^*(t)\\
       A(t) &B(t)
     \end{array}\right),
$$
where $A(t)$, $B(t)$, and $C(t)$ are $n\times n$ complex-valued
matrices, $B(t)$ and $C(t)$ are Hermitian matrices, and $A^*(t)$ is
the complex conjugate transpose of $A(t)$.
Then system $(1.1_{\lambda})$ can be written as
$$\begin{array}{rrll}
&&\Delta y_1(t)=A(t)y_1(t+1)+(B(t)+\lambda W_2(t))y_2(t),\nonumber\\
&&\Delta y_2(t)=(C(t)-\lambda W_1(t))y_1(t+1)-A^*(t)y_2(t),\quad t\in \mathcal{I}.
\end{array}                                                                                                               \eqno (1.2)$$

To ensure the existence and uniqueness of the solution of any
initial value problem for $(1.1_{\lambda})$, we always assume that
\begin{itemize}
\item[$(\mathbf{A}_1)$] $I_n-A(t)$ is invertible in $\mathcal{I}$.
\end{itemize}

It is known that $(1.1_{\lambda})$ contains the following formally
self-adjoint vector difference equation of order $2m$:
$$\sum_{j=0}^{m}(-1)^j\Delta^j[p_j(t)\Delta^jz(t-j)]=\lambda w(t)z(t),\quad t\in \mathcal{I},                    \eqno (1.3)$$
where $w(t)$ and $p_j(t), 0\le j\le m,$ are $l\times l$ Hermitian
matrices, $w(t)\geq 0$,  and $p_m(t)$ is invertible in $\mathcal{I}$.
The reader is referred to [28] for the details.

Spectral problems can be divided into two classifications.
Those defined over finite closed intervals with well-behaved coefficients are called regular;
otherwise they are called singular.

With the development of information technology and the wide applications of digital compute,
more and more discrete systems have appeared and they have attracted a lot of attention.
The study of fundamental theory of regular difference equations has a long history and
their spectral theory has formed a relatively complete theoretical system such as eigenvalue
problems, orthogonality of eigenfunctions and expansion theory (cf., [2, 17, 27, 36, 39, 41]).
Spectral problems for singular difference equations were firstly studied by Atkinson [2] in 1964,
and some significant progresses have been made since then (cf., e.g., [5, 6, 8, 16, 21, 22, 24, 25, 28, 32, 33, 37, 38]).
Especially, research on spectral theory of singular discrete Hamiltonian systems
has attracted a great deal of interest and some good results have been obtained (cf., [22, 24, 25, 28, 37, 38],
and references cited therein).
In 2006, the second author of the present paper established
the Weyl-Titchmarsh theory for system (1.1) with $a=0$ in [28].
Later, she with Ren studied the defect indices and definiteness conditions
and gave out complete characterizations of self-adjoint extensions for system (1.1) [24, 25].
Recently, she with Sun studied some spectral properties of system (1.1) [37].
These results have laid a foundation of our present research.

It is well known that regular discrete spectral problems have finite and then discrete spectra.
In particular, they can be transformed into eigenvalue problems of a special kind of matrices.
So they can be easily calculated by computer. Compared with regular problems, the spectral
set of a singular discrete spectral problem may contain some essential spectral points except
for isolated spectral points. Thus, it is difficult to study them. It is interesting to ask whether
the spectra of a singular spectral problem can be approximated by those of regular spectral
problems, and how to do it. Obviously, the study of regular approximations of spectra of
singular spectral problems plays an important role in both theory and practical applications.

Regular approximations of spectra of singular differential
equations have been investigated widely and deeply,
and some good results have been obtained, including spectral inclusion and spectral exactness [3, 4, 7, 18, 34, 35, 40, 43, 44].

To the best of our knowledge, there seem a few results about regular approximations of spectra of singular difference equations.
Recently, we studied this problem for singular second-order symmetric linear difference equations [19, 20].
For each self-adjoint subspace extension of a given singular second-order symmetric linear difference equation,
we constructed a sequence of regular problems and showed that the spectrum of the singular problem can be approximated
by the eigenvalues of this sequence.
Motivated by the ideas and methods used in [19, 20],
we shall study similar problems for singular discrete Hamiltonian system $(1.1_\lambda)$ in the present paper.
Although the methods are similar to that used in [19, 20], the problems investigated in the present paper are more complicated and difficult.
This results from the higher dimension and the partial shift operator $R$ in system $(1.1_\lambda)$.
We shall point out that there is another difficulty that will not be encountered in the continuous case.
It is that the maximal operator generated by (1.1) may be multi-valued,
and the corresponding minimal operator may be multi-valued or non-densely defined (see the detailed discussions in [24, 25, 29, 32]).
These facts were ignored in some existing literature including [28].
This is an essential difficulty that one would encounter in the study of the regular approximations of
spectra for difference expressions because the corresponding theory
of linear operators is not applicable in this case.

Fortunately, this major difficulty can be overcome by using the theory of linear subspaces (i.e., linear relations).
In 1961, Arens [1] initiated the study of linear relations, and his work was followed by many scholars [9-15].
Recently, some fundamental results of Hermitian subspaces including the Glazman-Krein-Naimark theory,
fundamental spectral properties of self-adjoint subspaces, and the resolvent convergence and spectral approximations of sequences of self-adjoint
subspaces were established [29-31]. A linear relation is actually a subspace in a
related product space, and obviously includes multi-valued and non-densely defined
linear operators in the related space.
Therefore, we shall study the regular approximations of spectra of system (1.1) in the framework of subspaces in a product space.

The rest of this paper is organized as follows.
In Section 2, some basic concepts and fundamental results about subspaces and system (1.1) are introduced,
including the maximal and minimal subspaces for (1.1), spectral inclusion, and spectral exactness.
In Section 3, the induced regular SSEs for any given SSE are constructed.
Section 4 pays attention to how to extend a subspace in the product space of the fundamental spaces on a proper subinterval
to a subspace in that on the original interval, i.e., how to do the ``zero extensions".
This problem can be very easily solved in the continuous case but hard in the discrete case.
Further, the invariance of spectral properties of the extended subspaces is given.
As a consequence, the extension from the induced regular
SSE to a subspace in the product space of the original Hilbert spaces is given,
and the invariance of spectral properties of the extended subspaces is obtained.
Regular approximations of spectra of system (1.1) in the limit circle case are studied in Section 5.
It is shown that the sequence of induced regular SSEs constructed in Section 3 is spectrally exact for any given SSE in this case.
In addition, it is obtained that the $k$-th eigenvalue of any given SSE
is exactly the limit of the $k$-th eigenvalues of the induced regular SSEs in this case.
Furthermore, error estimates for the approximations of eigenvalues are given in this case.
Section 6 is concerned with regular approximations of spectra of system (1.1) in the limit point and intermediate cases.
It is only shown that spectral inclusion holds in each case.
\medskip

\noindent{\bf Remark 1.1.}
We shall further study regular approximations of spectra of singular discrete linear Hamiltonian systems
with two singular endpoints in our forthcoming paper.
\medskip

\section{Preliminaries}

This section is divided into three parts. In Section 2.1, we recall some basic concepts and
fundamental results about subspaces.
In Section 2.2, we first introduce the maximal, pre-minimal, and minimal subspaces corresponding
to (1.1). Then, we list some useful results about (1.1), which will be used in the sequent sections.
Some useful results about resolvent convergence of sequences of self-adjoint subspaces are introduced in Section 2.3.
\medskip

\subsection{Some basic concepts and fundamental results about subspaces}

By ${\mathbf C}, {\mathbf R}$ and ${\mathbf Z^+}$ denote the sets of the complex numbers,
real numbers, and positive integer numbers, respectively.
Let $X$ be a complex Hilbert space with inner product $\langle
\cdot,\cdot\rangle$, and $T$ a linear subspace (briefly,
subspace) in the product space $X^2$ with the following induced inner product,
still denoted by $\langle \cdot,\cdot\rangle$ without any confusion:\vspace{-0.2cm}
$$\langle(x,f), (y,g)\rangle=\langle x, y\rangle+\langle f, g\rangle,
\;\;(x,f), (y,g)\in X^2.\vspace{-0.2cm}$$
The domain $D(T)$, range $R(T)$, and null space $N(T)$ of $T$ are respectively defined by
\vspace{-0.2cm}
$$\begin{array}{rrll}
D(T):&=&\{x\in X:\, (x,f)\in T \;{\rm for\; some}\;f\in X\},\\[0.4ex]
R(T):&=&\{f\in X:\, (x,f)\in T \;{\rm for\; some}\;x\in X\},\\[0.4ex]
N(T):&=&\{x\in X:\, (x,0)\in T \}.
\end{array}\vspace{-0.2cm}$$
Its adjoint subspace $T^*$ is defined by\vspace{-0.2cm}
$$T^*=\{(y,g)\in X^2:\,\langle f,y\rangle=\langle x,g\rangle\;{\rm
for\;all}\; (x,f)\in T\}.\vspace{-0.2cm}$$
Further, denote\vspace{-0.2cm}
$$T(x):=\{f\in X:\,(x,f)\in T\},\;\;T^{-1}:=\{(f,x):\,(x,f)\in T\}.\vspace{-0.1cm}$$

It is evident that $T(0)=\{0\}$ if and only if $T$ can uniquely determine a
single-valued linear operator from $D(T)$ into $X$ whose graph is $T$.
A single-valued linear operator is briefly called a linear operator.
For convenience, a linear operator in $X$ will always be identified with
a subspace in $X^2$ via its graph.

A subspace $T\subset X^2$ is called a Hermitian subspace if $T\subset T^*$, and it is called a self-adjoint
subspace if $T=T^*$. A Hermitian subspace $S$ is called a Hermitian subspace extension of $T$ if $T\subset S$,
and it is called a self-adjoint subspace extension of $T$ if $T\subset S$ and $S$ is a self-adjoint subspace.
In addition, a subspace $T$ is a Hermitian subspace if and only if
$\langle f, y\rangle=\langle x, g\rangle$ for all $(x,f), (y,g)\in T$.

Let $T$ and $S$ be two subspaces in $X^2$ and $\al \in {\mathbf C}$.
Define\vspace{-0.2cm}
$$\begin{array}{cl}
\alpha T:=\{(x,\alpha f):\, (x,f)\in T\},\\[0.6ex]
T+S:= \{(x,f+g):\,(x,f)\in T, (x,g)\in S\},\\[0.6ex]
ST:=\{(x,g)\in X^2:\, (x,f)\in T, (f,g)\in S \;{\rm for\; some}\;f\in X\}.
\end{array}\vspace{-0.2cm}$$
It is evident that if $T$ is closed, then $T-\lambda I_{id}$ is closed and $(T-\lambda I_{id})^*=T^*-{\bar \lambda}I_{id}$,
where $I_{id}:=\{(x,x):\, x\in X\}$, without any confusion we briefly denote it by $I$.
\medskip

For the following definition, the reader is referred to [15, 30, 31].
\medskip

\noindent{\bf Definition 2.1.}   Let $T$ be a subspace in $X^2$.
\begin{itemize}
\item[{\rm (1)}]The set $\rho(T):=\{\ld\in {\mathbf C}:\,(\ld I-T)^{-1}\;{\rm is\;  a\; bounded\;
 linear\; operator\; defined\; on}\;X\}$
 is called the resolvent set of $T$.

\item[{\rm (2)}]The set $\si(T):={\mathbf
C}\setminus \rho(T)$ is called the spectrum of $T$.
\end{itemize}
\medskip

\noindent\textit{{\bf Lemma 2.1 {\rm[30, Lemma 2.1]}.} Let $T$ be a closed subspace in $X^2$.
Then
$$\begin{array}{rrll}
\rho(T^{-1})\setminus \{0\}&=&\{\lambda^{-1}:\;\lambda\in \rho(T)\;{\rm with}\;\lambda\ne 0\},\\[0.4ex]
\sigma(T^{-1})\setminus \{0\}&=&\{\lambda^{-1}:\;\lambda\in \sigma(T)\;{\rm with}\;\lambda\ne 0\}.
\end{array}                                                                                  $$
Consequently, if $\rho(T)\ne \emptyset$, then
$$\sigma((\lambda_0 I-T)^{-1})\setminus \{0\}=\{(\lambda_0-\lambda)^{-1}:\;\lambda\in \sigma(T)\},\;\;
\lambda_0\in \rho(T).                    $$}
\vspace{-0.6cm}

\noindent\textit{{\bf Lemma 2.2 {\rm [19, Lemma 2.1]}.} Let $T$ be a closed subspace in $X^2$.
Then $\lambda\in \rho(T)$ if and only if $R(\lambda I-T)=X$ and $N(\lambda I-T)=\{0\}.$}
\medskip

\noindent\textit{{\bf Lemma 2.3 {\rm [30, Theorem 3.6]}.} Assume that $X_1$ is a proper closed subspace in $X$,
$P:\, X\to X_1$ the orthogonal
projection, and $T$ a self-adjoint subspace in $X_1^2$. Then
\begin{itemize}\vspace{-0.2cm}
\item[{\rm (i)}] $T'=TG(P)$ is a self-adjoint subspace in $X^2$
with $D(T')=D(T)\oplus X_1^\perp$;\vspace{-0.2cm}
\item[{\rm (ii)}] $\si(T')=\si(T)\cup\{0\}$.
\end{itemize}}

\subsection{Maximal, pre-minimal, and minimal subspaces}

In this subsection, we first introduce the concepts of maximal, pre-minimal, and minimal subspaces, and then list
some useful results about system $(1.1_{\lambda})$.

For any integer interval $\mathcal{\mathcal{I}}=\{t\}_{t=a}^{b}$ with $-\infty<a<b\leq +\infty$, we denote
$$\mathcal{I}^+:=\{t\}_{t=a}^{b+1},\;\;\;\; l(\mathcal{I}):=\{y:y=\{y(t)\}_{\mathcal{I}^+}\subset \mathbf{C}^{2n}\},$$
where $b+1$ means $+\infty$ in the case of $b=+\infty$.
Denote\vspace{-0.1cm}
$$
\mathcal {L}^2_{W}(\mathcal{I}):=\left\{y\in l(\mathcal{I}):\sum_{t\in \mathcal{I}}R(y)^*(t)W(t)R(y)(t)<+\infty\right\}
\vspace{-0.1cm}$$
with the semi-scalar product\vspace{-0.1cm}
$$
\langle x,y\rangle:=\sum_{t\in \mathcal{I}}R^*(y)(t)W(t)R(x)(t).
\vspace{-0.1cm}$$
Further, we define $\|y\|:=(\langle y,y\rangle)^{1/2}$ for $y\in
\mathcal {L}^2_{W}(\mathcal{I})$. Since the weighted function $W(t)$ may be
singular in $\mathcal{I}$, $\|\cdot\|$ is a semi-norm. We denote\vspace{-0.1cm}
$$
L^2_{W}(\mathcal{I}):=\mathcal {L}^2_W(\mathcal{I})/\{y\in\mathcal {L}^2_W(\mathcal{I}): \|y\|=0\}.
\vspace{-0.1cm}$$
Then $L^2_W(\mathcal{I})$ is a Hilbert space with the inner product $\langle\cdot,\cdot\rangle$ (cf. [28, Lemma 2.5]).
For a function $y\in \mathcal {L}^2_W(\mathcal{I})$, we denote by
$\tilde{y}$ the corresponding equivalent class in $L^2_W(\mathcal{I})$.
And for any $\tilde{y}\in L^2_W(\mathcal{I})$, by $y\in \mathcal {L}^2_W(\mathcal{I})$ denote
a representative of $\tilde{y}$. It is evident that
$\langle \tilde{x}, \tilde{y}\rangle=\langle x, y\rangle$ for any $\tilde{x}, \tilde{y}\in L^2_W(\mathcal{I})$.
Set\vspace{-0.1cm}
$$\begin{array}{rrll}
&&\mathcal {L}^2_{W,0}(\mathcal{I}):=\{y\in \mathcal {L}^2_{W}(\mathcal{I}):{\rm there
\; exist \;two \;integer}\; s,k\in \mathcal{I} \;{\rm
with}\;{s\le k}\nonumber\\
&& \quad \quad\quad \quad \quad \quad{\rm such\; that}\; y(t)=0
\;{\rm for}\; t\leq s\; {\rm and}\; t\ge k+1\}.
\end{array}\vspace{-0.1cm}$$
The natural difference operator corresponding to system $(1.1_{\lambda})$ is\vspace{-0.1cm}
$$
\mathscr{L}(y)(t):=J\Delta y(t)-P(t)R(y)(t).\vspace{-0.1cm}$$
Set\vspace{-0.1cm}
$$\begin{array}{rrll}
&&H:=\{(\tilde{y},\tilde{g})\in (L^2_W(\mathcal{I}))^2: {\rm there
\;exists\;}
y\in\tilde{y}\; {\rm such \;that}\nonumber\\
&&\quad \quad \quad \quad \quad \mathscr{L}(y)(t)=W(t)R(g)(t),\; t\in \mathcal{I}\}, \\[1.5ex]
&&H_{00}:=\{(\tilde{y},\tilde{g})\in H: {\rm there \; exists\;} y\in
\tilde{y}\;{\rm such\; that}\;  y\in\mathcal
{L}^2_{W,0}(\mathcal{I})\nonumber\\
&&\quad\quad\quad\quad{\rm and}\; \mathscr{L}(y)(t)=W(t)R(g)(t),\; t\in \mathcal{I}\},
\end{array}\vspace{-0.1cm}$$
$$H_0:=\overline{H}_{00},\vspace{-0.1cm}$$
where $H, H_{00}$, and $H_{0}$ are called the maximal, pre-minimal, and minimal subspaces
corresponding to system (1.1), respectively. By [24, Theorem 3.1], $H_{00}^*=H_{0}^*=H$,
which implies that $H_{0}$ is a closed Hermitian subspace in $(L^2_W(\mathcal{I}))^2$.
\medskip

By $n_{\lambda}$  denote the number of linearly independent square
summable solutions of $(1.1_\lambda)$ in $\mathcal{L}^2_W(\mathcal{I})$,
and by $d_{\lambda}$ denote the defect index of $H_0$ and $\bar{\lambda}$.
By [24, Corollary 5.1] we know that $ n_{\lambda}=d_{\lambda}$ if and only if the following definiteness condition is satisfied:
\begin{itemize}\vspace{-0.2cm}
\item [\rm {$(\mathbf{A}_2)$}]  There exists a finite subset $\mathcal{I}_0:=[s_0, t_0]\subset \mathcal{I}$ such that for some
$\lambda\in\mathbf{C}$, any non-trivial solution $y(t)$ of $(1.1_{\lambda})$
satisfies
\begin{eqnarray*}
\sum_{t\in \mathcal{I}_0}R(y)^*(t)W(t)R(y)(t)>0.
\end{eqnarray*}\vspace{-0.2cm}
\end{itemize}
\vspace{-0.2cm}

\noindent{\bf Remark 2.1.}
\begin{itemize}
\item [\rm {(1)}] It has been shown in [24] that
if the inequality in $(\mathbf{A}_2)$ holds for some $\lambda\in \mathbf C$,
then it holds for all $\lambda\in \mathbf{C}$.
Several sufficient conditions for $(\mathbf{A}_2)$ were given in [24].
Furthermore, it was pointed out that $H_0$ may be non-densely defined or multi-valued
in [25, Section 6].

\item [\rm {(2)}] It has been shown that $(\mathbf{A}_2)$ is equivalent to that
for any $(\tilde{y},\tilde{g})\in H$, there exists a
unique $y\in \tilde{y}$ such that $\mathscr{L}(y)(t)=W(t)R(g)(t)$ for $t\in \mathcal{I}$ in [24, Theorem 4.2].
In this case, we briefly write $(y,\tilde{g})\in H$ in the rest of the paper.

\item [\rm {(3)}] Even if $(\mathbf{A}_1)$ and $(\mathbf{A}_2)$ hold, $L^2_W(\mathcal{I})$ may be finite-dimensional
since $W(t)\geq 0$ for $t\in \mathcal{I}$.
In this case, $d_{\lambda}\equiv 2n$ for any $\lambda \in\mathbf{C}$.
\end{itemize}
\medskip

In the sequel, it is always assumed that $(\mathbf{A}_2)$ holds. It has been shown by [28, Corollary
4.1] that $n_\lambda\geq n$ for each $\lambda\in \mathbf{C}\backslash \mathbf{R}$.
Let $d_{\pm}$ be the positive and negative indices of $H_0$.
Since $n_\lambda\leq 2n$ and $n_\lambda=d_\lambda$ for each $\lambda\in \mathbf{C}$, we have
$n\leq d_{\pm}\leq 2n$.
By [9, Corollary of Theorem 15 and Theorem 18], $H_0$ has an SSE in $L^2_W(\mathcal{I})$ if and only if $d_+= d_-$.
So we always assume that the following holds in the sequel:
\begin{itemize}
\item[$(\mathbf{A}_3)$] $d_+ =d_-=:d$.
\end{itemize}
In the minimal deficiency case of $d=n$, $\mathscr{L}$ is said to be in the limit point case (l.p.c.) at $t=+\infty$
and in the maximal deficiency case of $d=2n$, $\mathscr{L}$ is said to be in the limit circle case (l.c.c.) at $t=+\infty$.
We refer to the cases when $n<d<2n$ as $\mathscr{L}$ in the intermediate cases.
\medskip

Next, for any $x,y\in l(\mathcal{I})$, we denote\vspace{-0.1cm}
$$
(x,y)(t)=y^*(t)J x(t).
\vspace{-0.1cm}$$
In the case of $b=+\infty$, if $\lim_{t\to b}(x,y)(t)$ exists and is finite,
then its limit is denoted by $(x,y)(+\infty)$.

By [28, Lemma 2.1], one has that for any $x,y\in l(\mathcal{I})$ and any $s,k\in \mathcal{I}$,\vspace{-0.2cm}
$$
\sum_{t=s}^{k}[R(y)^*(t)\mathscr{L}(x)(t)- \mathscr{L}(y)^*(t)R(x)(t)]= (x,y)(t)|_{s}^{k+1}. \vspace{-0.2cm}\eqno(2.1)$$
Hence, for $(x,\tilde{f})$, $(y,\tilde{g})\in H$, we get from (2.1) that
$$\begin{array}{rrll}
&&\sum_{t=s}^{k}[R(y)^*(t)W(t)R(f)(t)- R(g)^*(t)W(t)R(x)(t)]\\[1.0ex]
&=&\sum_{t=s}^{k}[R(y)^*(t)\mathscr{L}(x)(t)- \mathscr{L}(y)^*(t)R(x)(t)]\\[1.0ex]
&=&(x,y)(t)|_{s}^{k+1},
\end{array}$$
which yields that $\lim_{t\to +\infty}(x,y)(t)$ exists and is finite for all $(x,\tilde{f})$, $(y,\tilde{g})\in H$.
Further, by [28, Theorem 2.1] we get that for any $\lambda\in \mathbf{C}$, $c_0\in \mathcal{I}$, and any
solutions $y_{\lambda}(t)$ and $y_{\bar{\lambda}}(t)$ of $(1.1_{\lambda})$ and $(1.1_{\bar{\lambda}})$, respectively,\vspace{-0.1cm}
$$
(y_{\lambda}, y_{\bar{\lambda}})(t)=(y_{\lambda}, y_{\bar{\lambda}})(c_0), \quad t\in \mathcal{I}^+.  \vspace{-0.1cm}  \eqno(2.2)$$

\noindent\textit{{\bf Lemma 2.4 {\rm [25, Lemma 3.3]}.} Assume that $(\mathbf{A}_1)$ and $(\mathbf{A}_2)$ hold.
Then for any given finite subset $\tilde{\mathcal{I}}=\{t\}_{t=s}^{k}$ with $\mathcal{I}_0\subset \tilde{\mathcal{I}}\subset \mathcal{I}$
and for any given $\alpha, \beta\in \mathbf{C}^{2n}$, there exists $g=\{g(t)\}_{t=s}^{k+1}\in \mathbf{C}^{2n}$
such that the following boundary value problem:
$$\begin{array}{rrll}
&&\mathscr{L}(y)(t)=W(t)R(g)(t),\,t\in \tilde{\mathcal{I}},\\
&&y(s)=\alpha,\qquad y(k+1)=\beta,\\
\end{array}$$
has a solution $y=\{y(t)\}_{t=s}^{k+1}\in \mathbf{C}^{2n}$.
}\medskip

The following four lemmas are about SSE of $H_0$ and will be used
in constructing proper induced regular SSEs for any
given SSE of $H_{0}$.\medskip

\noindent\textit{{\bf Lemma 2.5 {\rm [25, Theorem 5.12]}.} Assume that $(\mathbf{A}_1)$ and $(\mathbf{A}_2)$
hold and $\mathcal{I}=\{t\}_{t=a}^b$ is finite. Then a subspace $H_1\subset (L^2_W(\mathcal{I}))^2$ is
an SSE of $H_{0}$ if and only if there exist two $2n\times 2n$ matrices $M$ and $N$ such that\vspace{-0.1cm}
$${\rm rank}\, (M,N)=2n,\qquad MJM^*=NJN^*,                                               \vspace{-0.1cm}$$
$$\vspace{-0.1cm}
H_1=\{(y,\tilde{g})\in H: My(a)-Ny(b+1)=0\}.
                                                                                            \vspace{-0.1cm}\eqno (2.3)$$
}

\noindent\textit{{\bf Lemma 2.6 {\rm [25, Theorem 5.10]}.} Assume that $(\mathbf{A}_1)$, $(\mathbf{A}_2)$, and $(\mathbf{A}_3)$ hold
and $\mathscr{L}$ is in l.c.c. at $t=+\infty$.
Let $\theta_1,\theta_2,\ldots,\theta_{2n}$ be $2n$
linearly independent  solutions of $(1.1_{\lambda})$ with $\lambda\in \mathbf R$ and satisfy
the following initial condition:\vspace{-0.1cm}
$$(\theta_1,\theta_2,\ldots,\theta_{2n})(a,\lambda)=I_{2n}.                               \vspace{-0.1cm}\eqno (2.4)$$
Then a subspace $H_1\subset (L^2_W(\mathcal{I}))^2$ is an SSE of $H_{0}$ if and
only if there exist two $2n\times 2n$ matrices $M$ and $N$ such that\vspace{-0.1cm}
$${\rm rank}\, (M,N)=2n,\qquad MJM^*=NJN^*,                                                  \vspace{-0.1cm}\eqno (2.5)$$
$$\vspace{-0.1cm}
H_1=\{(y,\tilde{g})\in H:
My(a)-N\left(\begin{array}{c}
                     (y,\theta_1)(+\infty) \\
                     \vdots \\
                     (y,\theta_{2n})(+\infty)
                   \end{array}
\right)=0\}.
                                                                                           \vspace{-0.1cm}\eqno (2.6)$$
}

\noindent\textit{{\bf Lemma 2.7 {\rm [25, Theorem 5.9]}.} Assume that $(\mathbf{A}_1)$, $(\mathbf{A}_2)$, and $(\mathbf{A}_3)$ hold
and $\mathscr{L}$ is in l.p.c. at $t=+\infty$. Then a subspace $H_1\subset (L^2_W(\mathcal{I}))^2$ is an SSE of $H_{0}$ if and
only if there exists a matrix $M_{n\times 2n}$ satisfying the self-adjoint conditions:\vspace{-0.2cm}
$$
{\rm rank}\, M=n,\quad MJM^*=0                                                              \vspace{-0.2cm}\eqno (2.7)$$
such that $H_1$ can be defined by\vspace{-0.2cm}
$$
H_1=\{(y,\tilde{g})\in H: My(a)=0\}.                                                         \vspace{-0.2cm}\eqno (2.8)$$}
\vspace{-0.4cm}

\noindent\textit{{\bf Lemma 2.8 {\rm [25, Theorem 5.8]}.} Assume that $(\mathbf{A}_1)$, $(\mathbf{A}_2)$, and $(\mathbf{A}_3)$ hold
and $\mathscr{L}$ is in the intermediate case at $t=+\infty$;
that is, $n<d<2n$. And assume that there exists $\lambda_0\in \mathbf{R}$ such that
system (1.1) has $d$ linear independent solutions $\psi_1, \cdots, \psi_d$ in $\mathcal{L}^2_W(\mathcal{I})$.
Let them be arranged such that\vspace{-0.1cm}
$$\Lambda:=((\psi_i, \psi_j)(+\infty))_{1\leq i,j\leq 2d-2n}=((\psi_i, \psi_j)(a))_{1\leq i,j\leq 2d-2n}\vspace{-0.1cm}$$
is invertible.
Then a subspace $H_1\subset (L^2_W(\mathcal{I}))^2$ is an SSE of $H_{0}$ if and
only if there exist two matrices $M_{d\times 2n}$ and $N_{d\times (2d-2n)}$ such that\vspace{-0.1cm}
$${\rm rank}\, (M,N)=d,\qquad MJM^*=N\Lambda^{\rm T} N^*,                                     \vspace{-0.1cm}\eqno (2.9)$$
and\vspace{-0.1cm}
$$
H_1=\{(y,\tilde{g})\in H:
My(a)-N\left(\begin{array}{c}
                     (y,\psi_1)(+\infty) \\
                     \vdots \\
                     (y,\psi_{2d-2n})(+\infty)
                   \end{array}
\right)=0\}.
                                                                                        \vspace{-0.1cm} \eqno (2.10)$$
}

\noindent{\bf Remark 2.2.}
By [25, Theorem 4.2] one can rearrange $\psi_1, \cdots, \psi_d$ such that $\Lambda$ is invertible,
where $\psi_1, \cdots, \psi_d$ and $\Lambda$ are specified in Lemma 2.8.
\medskip

\subsection{Resolvent convergence, spectral inclusion, and spectral exactness}

In this subsection, we recall some basic concepts,
including spectral inclusion, spectral exactness,
and strong resolvent convergence for self-adjoint subspaces and list some useful results.\medskip
\vspace{-0.3mm}

\noindent{\bf Definition 2.2 {\rm[30, Definition 4.1]}.}  Let $\{T_k\}_{k=1}^\infty$ and $T$ be self-adjoint subspaces in $X^2$.
$\{T_k\}_{k=1}^\infty$ is said to converge to $T$ in the strong resolvent sense
(briefly, SRC) if for some $\lambda\in  {\mathbf C}\setminus{\mathbf R}$, $(\lambda I-T_k)^{-1}$ is strongly convergent to $(\lambda I-{T})^{-1}$;
that is, $\|(\lambda I-T_k)^{-1}f-(\lambda I-T)^{-1}f\|\to 0$ as $k\to \infty$ for any $f\in X,$ denoted by
$(\lambda I-T_k)^{-1} \stackrel{s} {\to} (\lambda I-T)^{-1}$.
\medskip

\noindent{\bf Definition 2.3 {\rm[30, Definition 5.1]}.} Let $\{T_k\}_{k=1}^{\infty}$  and $T$ be subspaces in $X^2$.
\begin{itemize}
\item[{\rm (1)}] The sequence $\{T_k\}_{k=1}^{\infty}$ is said to be spectrally inclusive for
$T$ if for any $\lambda\in\sigma(T)$, there exists a sequence
$\{\lambda_k\}_{k=1}^{\infty},\;\lambda_k\in\sigma(T_k)$, such that $\lim_{k\to \infty}\lambda_k=\lambda$.
\item[{\rm (2)}] The sequence $\{T_k\}_{k=1}^{\infty}$ is said to be spectrally exact for $T$
if it is spectrally inclusive and every limit point of
any sequence $\{\lambda_k\}_{k=1}^{\infty}$ with $\lambda_k\in \sigma(T_k)$ belongs to $\sigma(T)$.
\end{itemize}

The following result gives a sufficient condition
for resolvent convergence of sequences of self-adjoint subspaces in the strong sense.\medskip

\noindent\textit{{\bf Lemma 2.9 {\rm [30, Theorem 4.2]}.}  Let $\{T_k\}_{k=1}^{\infty}$ and $T$ be self-adjoint subspaces in $X^2$.
Then $\{T_k\}$ is SRC to $T$ if $T$ has a core $T_0$ satisfying that
$T_0=\lim_{k\to \infty}T_k$; that is, for any $(x,f)\in T_0$, there exists $(x_k,f_k)\in T_k$
such that $(x,f)=\lim_{k\to\infty}(x_k,f_k).$}\medskip

A subspace $T_0$ is called a core of a closed subspace $T$ if $\overline{T}_0 = T$ (see
Definition 3.3 in [29]).\medskip

The following result gives a sufficient condition for spectral inclusion and
spectral exactness of a sequence of self-adjoint subspaces, which will take an important role in
the study of regular approximations of spectrum.\medskip

\noindent\textit{{\bf Lemma 2.10 {\rm [30, Theorem 5.4]}.}  Let $X_k, k\geq 1,$ be proper closed subspaces in $X$,
$P_k: X\to X_k$ orthogonal projections,
and $T$ and $\{T_k\}_{k=1}^{\infty}$ self-adjoint subspaces in $X^2$
and $X_k^2$, respectively.
Assume that $0\not\in \sigma(T)\ne \emptyset$ and
$\sigma(T_k)\ne \emptyset$ for $k\ge 1$, and set $T_k':=T_kG(P_k)$.
If $\{T_k'\}_{k=1}^\infty$ is SRC to $T$,
then $\{T_k\}_{k=1}^\infty$ is spectrally inclusive for $T$.
Further, if for any $\lambda \in \mathbf{C}\setminus {\mathbf R}$,
$\|(\lambda I-T_k)^{-1}G(P_k)-(\lambda I-T)^{-1}\| \to 0$ as $k \to 0$, denoted by
$(\lambda I-T_k)^{-1}G(P_k)\stackrel{n}{\to} (\lambda I-T)^{-1}$,
then $\{T_k\}_{k=1}^\infty$ is spectrally exact for $T$.}
\medskip

\section{Constructing induced regular self-adjoint subspace extensions}

Let $\mathcal{I}_r=\{t\}_{t=a}^{b_r},$ where $a<t_0<b_r<+\infty,
b_r\leq b_{r+1}, r\in \mathbf{Z^+}$, and
$b_r\rightarrow +\infty$ as $r\rightarrow \infty$,
where $t_0$ is specified by $(\mathbf{A}_2)$.
For convenience, by $H^r$ and $H^r_0$ denote
the corresponding maximal and minimal subspaces
corresponding to system (1.1) or $\mathscr{L}$
on $\mathcal{I}_r$, respectively.

Our main object in this section is to construct proper
induced regular SSEs $H_{1,r}$ of $\mathscr{L}$ on $\mathcal{I}_r$ for any given SSE $H_1$ of $H_0$.
We shall use the spectra of  $H_{1,r}$ to approximate the
spectrum of the given SSE $H_1$.
The discussions are divided into the following three cases:
$\mathscr{L}$ is in l.c.c., l.p.c., and the intermediate cases at $t=+\infty$.
\medskip

\noindent{\bf Case 1. The limit circle case}
\medskip

Let $\mathscr{L}$ be in l.c.c. at $t=+\infty$.
And let $\theta_1,\theta_2,\ldots,\theta_{2n}$ be defined in Lemma 2.6.
Set\vspace{-0.1cm}
$$\Theta(t,\lambda)=(\theta_1(t,\lambda),\theta_2(t,\lambda),\ldots,\theta_{2n}(t,\lambda)).          \vspace{-0.1cm}\eqno (3.1)$$
Then by (2.2) and (2.4) we get that\vspace{-0.1cm}
$$\Theta^*(t,\lambda)J\Theta(t,\lambda)= J.                                                                \vspace{-0.1cm}\eqno (3.2)$$

Suppose that $H_1$ is any fixed SSE
of $H_0$ and characterized by (2.6),
and matrices $M, N$ satisfy (2.5). Let\vspace{-0.1cm}
$$JM^*=(\rho_1,\rho_2,\ldots,\rho_{2n}),\quad N=(n_{ij})_{{2n}\times {2n}},                   \vspace{-0.1cm}\eqno (3.3) $$
and $\varphi_i:=\sum\limits^{2n}_{j=1}\bar{n}_{ij}\theta_j, 1\leq i\leq 2n$.
It is evident that $\varphi_i\in D(H), 1\leq i\leq 2n$.
By Lemma 2.4 there exist $\beta_i:=(\omega_i, \tilde{\tau}_i)\in H$ ($1\leq i\leq 2n$) such that\vspace{-0.1cm}
$$
\omega_i(a)=\rho_i,\quad  \omega_i(t)=\varphi_i(t),\quad t\geq t_0+1,                             \vspace{-0.1cm}\eqno (3.4)
$$
where $t_0$ is specified by $(\mathbf{A}_2)$.
By noting that\vspace{-0.2cm}
$$
My(a)=(JM^*)^*J y(a)=\left(\begin{array}{c}
                                   \omega_1^*(a) \\
                                   \vdots \\
                                   \omega_{2n}^*(a)
                                 \end{array}\right)Jy(a)=\left(\begin{array}{c}
                                   (y,\omega_1)(a) \\
                                   \vdots \\
                                   (y,\omega_{2n})(a)
                                 \end{array}\right),
$$
$$
N\left(\begin{array}{c}(y,\theta_1)(+\infty)\\\vdots\\
(y,\theta_{2n})(+\infty)\end{array}\right)
=\left(\begin{array}{c}(y,\sum_{j=1}^{2n}\bar{n}_{1j}\theta_j)(+\infty)\\\vdots\\
(y,\sum_{j=1}^{2n}\bar{n}_{(2n)j}\theta_j))(+\infty)\end{array}\right)=\left(\begin{array}{c}(y,\omega_1)(+\infty)\\\vdots\\
(y,\omega_{2n})(+\infty)\end{array}\right),
$$
$H_1$ in (2.6) can be rewritten as the following form:\vspace{-0.2cm}
$$H_1=\left\{(y,\tilde{g})\in H:\,\left(\begin{array}{c}(y,\omega_1)(a) \\ \vdots \\(y,\omega_{2n})(a)\end{array}\right)
-\left(\begin{array}{c}(y,\omega_1)(+\infty) \\ \vdots \\(y,\omega_{2n})(+\infty)\end{array}\right)=0
\right\}.                                                                                      \vspace{-0.2cm}   \eqno (3.5)$$
It can be easily verified that the set
$\{\beta_i\}_{i=1}^{2n}$ is a GKN-set
for $\{ H_0, H_0^*\}$. For the definition of a GKN-set
of Hermitian subspaces, the reader is referred to [29, Definition 4.1].

Next, we construct a proper induced regular SSE for $\mathscr{L}$ on $\mathcal{I}_r$ corresponding
to the given SSE $H_1$.

Let $b=b_r, P=M$, and $Q=N\Theta^*(b_r+1)J$ in Lemma 2.5.
Then $(\mathbf{A}_2)$ for (1.1) on $\mathcal{I}_r$ holds.
Since $\Theta$ and $J$ are invertible, one has that\vspace{-0.2cm}
$${\rm{rank}}(P,Q)={\rm{rank}}(M,N)=2n.\vspace{-0.2cm}$$
By (3.2) we have\vspace{-0.2cm}
$$QJQ^*=N\Theta^*(b_r+1)J\Theta(b_r+1)N^*=NJN^*,\;\;PJP^*=MJM^*.\vspace{-0.2cm}$$
Therefore,\vspace{-0.2cm}
$$PJP^*=QJQ^*.\vspace{-0.2cm}$$
In addition, because\vspace{-0.2cm}
$$\left( \begin{array}{c}
                (y,\theta_1)(t) \\
                \vdots\\
                (y,\theta_{2n})(t)
              \end{array}\right)=\Theta^*(t)J y(t),
\vspace{-0.2cm}$$
the subspace\vspace{-0.2cm}
$$H_{1,r}=\left\{(y,\tilde{g})\in H^r: \,My(a)-N \left( \begin{array}{c}
                (y,\theta_1)(b_r+1) \\
                \vdots\\
                (y,\theta_{2n})(b_r+1)
              \end{array}\right)=0\right\}             \vspace{-0.2cm}               \eqno (3.6)$$
is an SSE of $H^r_0$ by Lemma 2.5.
With a similar argument to that used in the above discussion
for (3.5), one can easily get that $H_{1,r}$ can be rewritten as\vspace{-0.2cm}
$$H_{1,r}=\left\{(y,\tilde{g})\in H^r:\,\left(\begin{array}{c}(y,\omega_1)(a) \\ \vdots \\(y,\omega_{2n})(a)\end{array}\right)
-\left(\begin{array}{c}(y,\omega_1)(b_r+1) \\ \vdots \\(y,\omega_{2n})(b_r+1)\end{array}\right)=0
\right\}.                                                                                      \vspace{-0.2cm}   \eqno (3.7)$$
We call $H_{1,r}$ an induced regular SSE of $H_1$ on $\mathcal{I}_r$.
Further, it can be easily verified that
$\{{\beta}_i|_{\mathcal{I}_r^+}:=\{({\omega}_i,\tilde{\tau}_i)(t)\}_{t=a}^{b_r+1}\}_{i=1}^{2n}$
is a GKN-set for $\{H^r_0,{H^r_0}^*\}$.
\medskip

\noindent{\bf Case 2. The limit point case}
\medskip

Let $\mathscr{L}$ be in l.p.c. at $t=+\infty$.
Suppose that $H_1$ is any fixed SSE of $H_0$ and characterized by (2.8),
and the matrix $M_{n\times2n}$ satisfies (2.7).
Let\vspace{-0.2cm}
$$JM^*=(\alpha_{1},\ldots,\alpha_{n}).\vspace{-0.2cm}$$
By Lemma 2.4, there exist $\beta_{i}=(\omega_{i},\tilde{\tau}_{i})\in H$, $1\leq i\leq n$, satisfying\vspace{-0.2cm}
$$\omega_{i}(a)=\alpha_{i},\;\;\omega_{i}(t)=0,\;\;t\geq t_{0}+1,               \vspace{-0.2cm}          \eqno (3.8)$$
where $t_0$ is specified by $(\mathbf{A}_2)$.
It follows that\vspace{-0.2cm}
$$
My(a)=(JM^*)^*J y(a)=\left(\begin{array}{c}
                                   \omega_1^*(a) \\
                                   \vdots \\
                                   \omega_{n}^*(a)
                                 \end{array}\right)Jy(a)=\left(\begin{array}{c}
                                   (y,\omega_1)(a) \\
                                   \vdots \\
                                   (y,\omega_{n})(a)
                                 \end{array}\right),
$$
which implies that $H_1$ in (2.8) can be written as\vspace{-0.2cm}
$$H_1=\{(y,\tilde{g})\in H:\,\left(\begin{array}{c}
                                   (y,\omega_1)(a) \\
                                   \vdots \\
                                   (y,\omega_{n})(a)
                                 \end{array}\right)=0\}.            \vspace{-0.2cm}                  \eqno (3.9)$$
It is obvious that $\{\beta_{i}\}_{i=1}^{n}$ is
a GKN-set for $\{H_0,H_0^*\}$.

Now, we construct a proper regular SSE $H_{1,r}$, which is induced by $H_1$ on $\mathcal{I}_r$.
Set\vspace{-0.2cm}
$$P=\left( \begin{array}{c}
                M_{n\times 2n} \\
                0_{n\times 2n}
\end{array}\!\right),\;\;\;Q=-\left( \!\begin{array}{c}
                0_{n\times 2n} \\
                N_{n\times 2n}
              \end{array}\!\right)\Theta^*(b_r+1,\lambda)J,\vspace{-0.2cm}$$
where $N=(n_{ij})_{n\times 2n}$
with ${\rm rank} N=n$ and $NJN^*=0$, is any fixed matrix,
$\Theta$ is defined by (3.1), and $\lambda\in\mathbf{R}$ is any fixed number.
It can be easily verified that\vspace{-0.2cm}
$${\rm{rank}}(P,Q)=2n,\;\;PJP^*=QJQ^*=0.\vspace{-0.2cm}$$
Further, it follows that\vspace{-0.2cm}
$$\begin{array}{rrll}
P y(a)-Q y(b_r+1)&=&\left(\begin{array}{c}M_{n\times 2n}\\0_{n\times 2n}\end{array}\right)y(a)+
\left(\begin{array}{c}0_{n\times 2n} \\ N_{n\times 2n}\end{array}\!\right)\Theta^*(b_r+1)Jy(b_r+1)            \\[2.5ex]
&=&\left(\begin{array}{c}M_{n\times 2n}\\0_{n\times 2n}\end{array}\right)y(a)+
\left(\begin{array}{c}0_{n\times 2n} \\ N_{n\times 2n}\end{array}\!\right)
\left(\begin{array}{c}(y,\theta_1)(b_r+1) \\\vdots \\(y,\theta_{2n})(b_r+1)\end{array}\right)            \\[2.5ex]
&=&\left(\begin{array}{c}
     My(a) \\
N\left(\begin{array}{c}(y,\theta_1)(b_r+1) \\\vdots \\(y,\theta_{2n})(b_r+1)\end{array}\right)
\end{array}\right).
\end{array}$$
Therefore, by Lemma 2.5 one has that\vspace{-0.2cm}
$$H_{1,r}=\left\{(y,\tilde{g})\in H^r:\,My(a)=0,\,\,
N\left(\begin{array}{c}(y,\theta_1)(b_r+1) \\\vdots \\(y,\theta_{2n})(b_r+1)\end{array}\right)=0 \right \}\vspace{-0.2cm}  \eqno (3.10)$$
is an SSE of $H^r_0$. Let\vspace{-0.2cm}
$$\varphi_i:=\sum\limits^{2n}_{j=1}\bar{n}_{ij}\theta_j, \,1\leq i\leq n.                       \vspace{-0.2cm}             \eqno (3.11)$$
Similarly to the discussion for (3.5),
$H_{1,r}$ can be rewritten as\vspace{-0.2cm}
$$H_{1,r}=\left\{(y,\tilde{g})\in H^r:\,\left(\begin{array}{c}
(y,\omega_1)(a) \\\vdots \\(y,\omega_{n})(a)\end{array}\right)=0,\,
\left( \begin{array}{c}
(y,\varphi_1)(b_r+1) \\\vdots\\(y,\varphi_{n})(b_r+1)\end{array}\right)=0\right\}.                        \vspace{-0.2cm}   \eqno (3.12)$$
We call $H_{1,r}$ an induced regular SSE of $H_1$ on $\mathcal{I}_r$.
\medskip

\noindent{\bf Case 3. The intermediate cases}
\medskip

Let $\mathscr{L}$ be in the intermediate case at $t=+\infty$ with $n<d<2n$.
In the case, we always assume that
\begin{itemize}\vspace{-0.2cm}
\item [\rm {$(\mathbf{A}_4)$}] There exists $\lambda_0\in \mathbf{R}$ such that $(1.1_{\lambda_0})$
has $d$ linear independent solutions in $\mathcal{L}^2_W(\mathcal{I})$.
\end{itemize}
\vspace{-0.2cm}
Then we assert that $(1.1_{\lambda_0})$
has $d$ linear independent solutions $\psi_1, \cdots, \psi_d$ in $\mathcal{L}^2_W(\mathcal{I})$ such that\vspace{-0.1cm}
$$\begin{array}{rrll}
&&\Lambda:=((\psi_i, \psi_j)(+\infty))_{1\leq i,j\leq 2d-2n}\:{\rm is}\:{\rm a}\:{\rm diagonal}\:{\rm and}\:{\rm invertible}\:{\rm matrix};\\
&&((\psi_i, \psi_j)(+\infty))_{1\leq i,j\leq d}=\left(\begin{array}{cc}
                                                \Lambda                    & 0_{(2d-2n)\times (2n-d)} \\
                                                  0_{(2n-d)\times (2d-2n)} & 0_{(2n-d)\times (2n-d)}
                                              \end{array}\right).
\end{array}                                                                                      \vspace{-0.1cm}\eqno (3.13)$$
In fact, let $\tilde{\psi}_1, \cdots, \tilde{\psi}_d$ be any $d$ linear independent solutions of $(1.1_{\lambda_0})$
in $\mathcal{L}^2_W(\mathcal{I})$. Let $\tilde{\Psi}_1:=(\tilde{\psi}_1, \cdots, \tilde{\psi}_d)$. Then
$\tilde{\Psi}_1^*(t)J\tilde{\Psi}_1(t)=\tilde{\Psi}_1^*(+\infty)J\tilde{\Psi}_1(+\infty)$ by (2.2), which is a skew-Hermitian matrix.
In addition, ${\rm rank}(\tilde{\Psi}_1^*(+\infty)J\tilde{\Psi}_1(+\infty))=2d-2n$ by [25, Lemma 4.4]. Thus, there exists a unitary matrix $U$
such that\vspace{-0.1cm}
$$U^*(\tilde{\Psi}_1^*(+\infty)J\tilde{\Psi}_1(+\infty)) U=\left(\begin{array}{cc}
                                                \tilde{\Lambda}_{(2d-2n)\times (2d-2n)} & 0_{(2d-2n)\times (2n-d)} \\
                                                  0_{(2n-d)\times (2d-2n)}              & 0_{(2n-d)\times (2n-d)}
                                              \end{array}\right),\vspace{-0.1cm}$$
where $\tilde{\Lambda}_{(2d-2n)\times (2d-2n)}$ is a diagonal and invertible matrix.
Let $\Psi_1=(\psi_1, \cdots, \psi_d):=\tilde{\Psi}_1 U$. Then\vspace{-0.1cm}
$$((\psi_i, \psi_j)(+\infty))_{1\leq i,j\leq d}
=(\Psi_1^*(+\infty)J\Psi_1(+\infty))^{\rm T}
=\left(\begin{array}{cc}
                                                \tilde{\Lambda}_{(2d-2n)\times (2d-2n)} & 0_{(2d-2n)\times (2n-d)} \\
                                                  0_{(2n-d)\times (2d-2n)}              & 0_{(2n-d)\times (2n-d)}
                                              \end{array}\right)                                               \vspace{-0.1cm}$$
and so $\psi_1, \cdots, \psi_d$ are $d$ linear independent solutions of $(1.1_{\lambda_0})$ in $\mathcal{L}^2_W(\mathcal{I})$
and satisfy (3.13). Thus, this assertion holds.
In this case, we shall use these solutions $\psi_1, \cdots, \psi_d$ to characterize the self-adjoint subspace extensions $H_1$
of $H_0$ in Lemma 2.8.

Suppose that $H_1$ is any fixed SSE of $H_0$ and characterized by (2.10),
and matrices $M, N$ satisfy (2.9). Let\vspace{-0.1cm}
$$JM^*=(\gamma_1,\ldots,\gamma_{d}),\quad N=(n_{ij})_{{d}\times {(2d-2n)}},                         \vspace{-0.1cm}\eqno (3.14) $$
and set $\varphi_i:=\sum\limits^{2d-2n}_{j=1}\bar{n}_{ij}\psi_j, 1\leq i\leq d$.
Clearly, $\varphi_i\in D(H), 1\leq i\leq d$.
By Lemma 2.4 there exist $\beta_i:=(\omega_i, \tilde{\tau}_i)\in H$ ($1\leq i\leq d$) such that\vspace{-0.1cm}
$$
\omega_i(a)=\gamma_i,\quad  \omega_i(t)=\varphi_i(t),\quad t\geq t_0+1,                              \vspace{-0.1cm}\eqno (3.15)
$$
where $t_0$ is specified by $(\mathbf{A}_2)$.
Note that for any $y \in D(H)$, it follows that\vspace{-0.2cm}
$$
My(a)=(JM^*)^*J y(a)=\left(\begin{array}{c}
                                   \omega_1^*(a) \\
                                   \vdots \\
                                   \omega_{d}^*(a)
                                 \end{array}\right)Jy(a)=\left(\begin{array}{c}
                                   (y,\omega_1)(a) \\
                                   \vdots \\
                                   (y,\omega_{d})(a)
                                 \end{array}\right),
\vspace{-0.1cm}$$
$$\vspace{-0.1cm}
N\left(\begin{array}{c}
(y,\psi_1)(+\infty)\\
\vdots\\
(y,\psi_{2d-2n})(+\infty)\end{array}\right)
=\left(\begin{array}{c}
(y,\sum_{j=1}^{2d-2n}\bar{n}_{1j}\psi_j)(+\infty)\\
\vdots\\
(y,\sum_{j=1}^{2d-2n}\bar{n}_{dj}\psi_j))(+\infty)\end{array}\right)=\left(\begin{array}{c}
(y,\omega_1)(+\infty)\\\vdots\\
(y,\omega_{d})(+\infty)\end{array}\right).
$$
Hence, $H_1$ in (2.10) can be rewritten as the following form:\vspace{-0.2cm}
$$H_1=\left\{(y,\tilde{g})\in H:\,\left(\begin{array}{c}
(y,\omega_1)(a) \\
 \vdots \\
 (y,\omega_{d})(a)\end{array}\right)
-\left(\begin{array}{c}
(y,\omega_1)(+\infty) \\
\vdots \\
(y,\omega_{d})(+\infty)\end{array}\right)=0
\right\}.                                                                                     \vspace{-0.1cm}\eqno (3.16)$$
It can be easily verified that the set $\{\beta_i\}_{i=1}^{d}$ is a GKN-set for $\{ H_0, H_0^*\}$.

Next, we construct a proper induced regular SSE for $\mathscr{L}$ on $\mathcal{I}_r$ corresponding
to the given SSE $H_1$.

We still use the solutions $\psi_1, \cdots, \psi_d$ in $\mathcal{L}^2_W(\mathcal{I})$, which satisfy (3.13).
In addition, we add solutions $\psi_{d+1}, \cdots, \psi_{2n}$ such that
$\{\psi_{1}, \cdots, \psi_{2n}\}$ forms a basis of solutions of $(1.1_{\lambda_0})$.
Let $\Psi=(\psi_1, \cdots, \psi_{2n})$. Then $\Psi$ is obviously invertible. Set\vspace{-0.1cm}
$$P=\left( \begin{array}{c}
                M_{d\times 2n} \\
                0_{(2n-d)\times 2n}
\end{array}\!\right),\;\;\;Q=\left( \!
\begin{array}{ccc}
  N_{d\times (2d-2n)} & 0_{d\times (2n-d)} & 0_{d\times (2n-d)} \\
  0_{(2n-d)\times (2d-2n)} & I_{2n-d} & 0_{(2n-d)\times (2n-d)}
\end{array}
\!\right)\Psi^*(b_r+1)J.\vspace{-0.1cm}$$
It is obvious that\vspace{-0.1cm}
$${\rm{rank}}(P,Q)={\rm{rank}}(M,N)+2n-d=2n.\vspace{-0.1cm}$$
By (2.2), (3.13), and $MJM^*=N\Lambda^{\rm T} N^*$ we have\vspace{-0.1cm}
$$PJP^*=QJQ^*.\vspace{-0.1cm}$$
Further, it follows that\vspace{-0.1cm}
$$
P y(a)-Q y(b_r+1)
=\left(\begin{array}{c}
     My(a)-N\left(\begin{array}{c}
(y,\psi_1)(b_r+1)\\
\vdots\\
(y,\psi_{2d-2n})(b_r+1)\end{array}\right)\\
-\left(\begin{array}{c}(y,\psi_{2d-2n+1})(b_r+1) \\\vdots \\(y,\psi_{d})(b_r+1)\end{array}\right)
\end{array}\right).
\vspace{-0.1cm}$$
Therefore, by Lemma 2.5 one has that\vspace{-0.1cm}
$$H_{1,r}\!=\!\left\{\!\!(y,\tilde{g})\in H^r\!: My(a)\!-\!N\!\left(\!\!\begin{array}{c}
                                                    (y,\psi_1)(b_r+1)\\
                                                     \vdots\\
                                                    (y,\psi_{2d-2n})(b_r+1)\end{array}\!\!\right)\!\!=0,
\left(\!\!\begin{array}{c}(y,\psi_{2d-2n+1})(b_r+1) \\\vdots \\(y,\psi_{d})(b_r+1)\end{array}\!\!\right)\!\!=\!0 \!\right \}\vspace{-0.1cm}$$
is an SSE of $H^r_0$.
Similarly to the discussion for (3.5), one can easily get
that $H_{1,r}$ can be rewritten as\vspace{-0.2cm}
$$H_{1,r}\!=\!\left\{\!\!(y,\tilde{g})\in H^r\!:\!\left(\!\!\!\begin{array}{c}
(y,\omega_1)(a) \\\vdots \\(y,\omega_{d})(a)\end{array}\!\!\!\right)\!\!=\!\!\left(\!\!\!\begin{array}{c}
(y,\omega_1)(b_r+1) \\\vdots \\(y,\omega_{d})(b_r+1)\end{array}\!\!\!\right),
\left(\!\!\!\begin{array}{c}
(y,\psi_{2d-2n+1})(b_r+1) \\\vdots\\(y,\psi_{d})(b_r+1)\end{array}\!\!\!\right)\!\!\!=\!0 \!\right\},         \vspace{-0.2cm}   \eqno (3.17)$$
where $\omega_1,\ldots,\omega_{d}$ are defined by (3.15).
We call $H_{1,r}$ an induced regular SSE of $H_1$ on $\mathcal{I}_r$.
\medskip

\section{Extension of the induced regular self-adjoint subspace extensions to the whole space}

In this section, we first extend a subspace in the product space of the fundamental spaces on a proper subinterval
to a subspace in that on the original interval, and study spectral properties of the extended subspaces.
As a consequence, the extension from the induced regular SSE constructed in Section 3
to a subspace in $(L_W^2(\mathcal{I}))^2$ is given,
and the spectral properties of the extended subspaces are obtained.

Let $\mathcal{K}:=\{t\}_{t=a}^{b}$ be an integer interval,
where $a$ is a finite integer or $a=-\infty$ and $b$ is a finite integer or $b=+\infty$.
$\mathcal{K}^+$, $\mathcal{L}_W^2(\mathcal{K})$, and $L_W^2(\mathcal{K})$ can be well defined as in Section 2.2
with $\mathcal{I}$ replaced by $\mathcal{K}$.
For convenience, by $\langle \cdot,\cdot \rangle_{\mathcal{K}}$ and $\|\cdot\|_{\mathcal{K}}$
denote the inner product and norm of $(L_W^2(\mathcal{K}))^2$, respectively.

For any integer interval $\mathcal{J}\subsetneqq \mathcal{K}$, denote\vspace{-0.1cm}
$$\hat{L}_W^2(\mathcal{J}):=\{\tilde{y}\in L_W^2(\mathcal{K}):\;W(t)R(y)(t)=0,\;\;t\in \mathcal{K}\setminus \mathcal{J}\}.\vspace{-0.1cm}\eqno (4.1) $$

For any subspace $T$ in $(L_W^2(\mathcal{J}))^2$, denote\vspace{-0.1cm}
$$\begin{array}{rrll}
&&\hat{T}:=\{(\tilde{\hat{y}},\tilde{\hat{g}})\in (\hat{L}_W^2(\mathcal{J}))^2:\;
{\rm there\;exists}\;(\tilde{y},\tilde{g})\in T\;{\rm such\;that}\\[0.5ex]
&&\quad \quad \quad \quad \quad\quad \quad \quad \quad \quad\quad\quad
 \|y\|_{\mathcal{J}}=\|\hat{y}\|_{\mathcal{K}},\; \|g\|_{\mathcal{J}}=\|\hat{g}\|_{\mathcal{K}}\}.
 \end{array}                                                                                                   \vspace{-0.1cm}\eqno (4.2) $$

The following result can be easily verified, and so its details are omitted.\medskip

\noindent\textit{{\bf Proposition 4.1.}
$\hat{L}_W^2(\mathcal{J})$ is a closed subspace in $L_W^2(\mathcal{K})$ and
so $\hat{T}$ is a subspace in $(L_W^2(\mathcal{K}))^2$.
Moreover, $\langle x,y \rangle_{\mathcal{K}}=\langle x|_{\mathcal{J}^+},y|_{\mathcal{J}^+} \rangle_{\mathcal{J}}$
for any $\tilde{x}, \tilde{y}\in \hat{L}_W^2(\mathcal{J})$.}\medskip

\noindent\textit{{\bf Proposition 4.2.}
Let $T$ be a subspace in $(L_W^2(\mathcal{J}))^2$ and $\hat{T}$ be defined by (4.2). Then
\begin{itemize}\vspace{-0.2cm}
\item[{\rm (i)}] $T$ is a closed subspace in $(L_W^2(\mathcal{J}))^2$
                 if and only if $\hat{T}$ is a closed subspace in $(\hat{L}_W^2(\mathcal{J}))^2$;\vspace{-0.2cm}
\item[{\rm (ii)}] $T$ is a Hermitian subspace in $(L_W^2(\mathcal{J}))^2$
                 if and only if $\hat{T}$ is a Hermitian subspace in $(\hat{L}_W^2(\mathcal{J}))^2$;\vspace{-0.2cm}
\item[{\rm (iii)}] $T$ is a self-adjoint subspace in $(L_W^2(\mathcal{J}))^2$
                if and only if $\hat{T}$ is a self-adjoint subspace in $(\hat{L}_W^2(\mathcal{J}))^2$.
\end{itemize}}\medskip

\noindent{\bf Proof.}
(i) Assertion (i) can be directly derived from (4.2) and Proposition 4.1.

(ii) We first show the necessity.
Suppose that $T$ is a Hermitian subspace in $(L_W^2(\mathcal{J}))^2$.
For any $(\tilde{\hat{y}_i},\tilde{\hat{g}}_i)\in \hat{T},\,i=1,2$, by (4.2) there exist
$(\tilde{y}_i,\tilde{g}_i)\in T,\,i=1,2$, such that $\|y_i\|_{\mathcal{J}}=\|\hat{y}_i\|_{\mathcal{K}}$
and $\|g_i\|_{\mathcal{J}}=\|\hat{g}_i\|_{\mathcal{K}},\,i=1,2$, which are equivalent to
\vspace{-0.1cm}
$$W(t)R(y_i)(t)=W(t)R(\hat{y}_i)(t),\;\;W(t)R(g_i)(t)=W(t)R(\hat{g}_i)(t), \,t\in \mathcal{J},\,i=1,2.\vspace{-0.1cm}\eqno (4.3)$$
Since $T$ is Hermitian,\vspace{-0.1cm}
$$\langle g_1,y_2\rangle_{\mathcal{J}}=\langle y_1,g_2\rangle_{\mathcal{J}}.                           \vspace{-0.1cm}\eqno (4.4)$$
This, together with (4.1)-(4.3), yields that\vspace{-0.1cm}
$$\langle \hat{g}_1,\hat{y}_2\rangle_{\mathcal{K}}=\langle \hat{y}_1,\hat{g}_2\rangle_{\mathcal{K}}.   \vspace{-0.1cm}\eqno (4.5)$$
This implies that $\hat{T}$ is a Hermitian subspace in $(\hat{L}_W^2(\mathcal{J}))^2$.

Next we consider the sufficiency. Suppose that $\hat{T}$ is a Hermitian subspace in $(\hat{L}_W^2(\mathcal{J}))^2$.
For any $(\tilde{y}_i,\tilde{g}_i)\in T,\,i=1,2$, take $\tilde{\hat{y}}_i,\tilde{\hat{g}}_i\in \hat{L}_W^2(\mathcal{J}),\,i=1,2$,
such that (4.3) holds. Then $(\tilde{\hat{y}}_i,\tilde{\hat{g}}_i)\in \hat{T},\,i=1,2$.
Since $\hat{T}$ is Hermitian, one has that (4.5) holds.
This, together with (4.1) and (4.3), yields that (4.4) holds.
This implies that $T$ is a Hermitian subspace in $(L_W^2(\mathcal{J}))^2$.

(iii) We first show the necessity. Suppose that $T$ is a self-adjoint subspace in $(L_W^2(\mathcal{J}))^2$.
By (ii) we get that $\hat{T}$ is a Hermitian subspace in $(\hat{L}_W^2(\mathcal{J}))^2$.
So, it is only needed to show that $\hat{T}^*\subset\hat{T}$.
By the definition of adjoint subspace, for any given $(\tilde{\hat{y}_1},\tilde{\hat{g}}_1)\in \hat{T}^*$,
(4.5) holds for all $(\tilde{\hat{y}_2},\tilde{\hat{g}}_2)\in \hat{T}$.
Set $y_1:={\hat{y}_1}|_{\mathcal{J}^+}$ and $g_1:={\hat{g}_1}|_{\mathcal{J}^+}$.
Then $\tilde{y}_1,\tilde{g}_1\in L_W^2(\mathcal{J})$ and (4.3) holds for $i=1$.
In addition, for any $(\tilde{y}_2,\tilde{g}_2)\in T$, there exists
$(\tilde{\hat{y}_2},\tilde{\hat{g}}_2)\in \hat{T}$
such that (4.3) holds for $i=2$.
Hence, we get that $\langle \hat{g}_1,\hat{y}_2\rangle_{\mathcal{K}}=\langle g_1,y_2\rangle_{\mathcal{J}}$
and $\langle \hat{y}_1,\hat{g}_2\rangle_{\mathcal{K}}=\langle y_1,g_2\rangle_{\mathcal{J}}$.
It follows from (4.5) that\vspace{-0.1cm}
$$\langle g_1,y_2\rangle_{\mathcal{J}}=\langle y_1,g_2\rangle_{\mathcal{J}},\;\forall\, (\tilde{y}_2,\tilde{g}_2)\in T.\vspace{-0.1cm}\eqno (4.6)$$
Because $T$ is a self-adjoint subspace in $(L_W^2(\mathcal{J}))^2$, we get that $(\tilde{y}_1,\tilde{g}_1)\in T$.
Therefore, $(\tilde{\hat{y}_1},\tilde{\hat{g}}_1)\in \hat{T}$ by (4.2). This implies that $\hat{T}^*\subset\hat{T}$.
Thus, $\hat{T}$ is a self-adjoint subspace in $(\hat{L}_W^2(\mathcal{J}))^2$.

Next we consider the sufficiency. Suppose that $\hat{T}$ is a self-adjoint subspace in $(\hat{L}_W^2(\mathcal{J}))^2$.
Similarly, by (ii) we only need to show that $T^*\subset T$.
By the definition of adjoint subspace, for any given $(\tilde{y}_1,\tilde{g}_1)\in T^*$, (4.6) holds.
Take $\tilde{\hat{y}}_1,\tilde{\hat{g}}_1\in \hat{L}_W^2(\mathcal{J})$ such that (4.3) holds for $i=1$.
In addition, for any $(\tilde{\hat{y}_2},\tilde{\hat{g}}_2)\in \hat{T}$, by (4.2) there exists
$(\tilde{y}_2,\tilde{g}_2)\in T$ such that (4.3) holds for $i=2$.
It follows from (4.3) and (4.6) that (4.5) holds for all $(\tilde{\hat{y}_2},\tilde{\hat{g}}_2)\in \hat{T}$.
This implies that $(\tilde{\hat{y}_1},\tilde{\hat{g}}_1)\in \hat{T}^*=\hat{T}$.
This, together with (4.2) and (4.3) with $i=1$, yields that $(\tilde{y}_1,\tilde{g}_1)\in T$.
Therefore, $T^*\subset T$ and so $T$ is a self-adjoint subspace in $(L_W^2(\mathcal{J}))^2$.
The whole proof is complete.
\medskip

\noindent\textit{{\bf Proposition 4.3.}
Let $T$ be a closed subspace in $(L_W^2(\mathcal{J}))^2$ and $\hat{T}$ be defined by (4.2). Then
$\sigma(\hat{T})=\sigma(T)$.}\medskip

\noindent{\bf Proof.} By Lemma 2.2, it suffices to show that\vspace{-0.1cm}
$$N(\lambda I-\hat{T})\neq \{0\} \Leftrightarrow N(\lambda I-T)\neq \{0\}, \vspace{-0.1cm}                         \eqno(4.7)$$
$$R(\lambda I-\hat{T})=\hat{L}_W^2(\mathcal{J}) \Leftrightarrow R(\lambda I-T)=L_W^2(\mathcal{J}). \vspace{-0.1cm} \eqno(4.8)$$

We first show that (4.7) holds. Suppose that $N(\lambda I-T)\neq \{0\}$. Then
there exists $\tilde{x}\in D(T)$ with $\|x\|_{\mathcal{J}}>0$ such that
$(\tilde{x},0)\in \lambda I-T$, i.e., $(\tilde{x},\lambda \tilde{x})\in T$.
Take $\tilde{\hat{x}}\in \hat{L}_W^2(\mathcal{J})$ satisfying $\|\hat{x}\|_{\mathcal{K}}=\|x\|_{\mathcal{J}}$.
Then, by (4.2), $(\tilde{\hat{x}},\lambda \tilde{\hat{x}})\in \hat{T}$, i.e., $(\tilde{\hat{x}},0)\in \lambda I-\hat{T}$.
So $N(\lambda I-\hat{T})\neq \{0\}$ by $\|\hat{x}\|_{\mathcal{K}}>0$.
On the other hand, suppose that $N(\lambda I-\hat{T})\neq \{0\}$.
Then, there exists $\tilde{\hat{y}}\in D(\hat{T})$ with $\|\hat{y}\|_{\mathcal{K}}>0$ such that
$(\tilde{\hat{y}},0)\in \lambda I-\hat{T}$, which implies that
$(\tilde{\hat{y}},\lambda \tilde{\hat{y}})\in \hat{T}$.
By (4.2), there exists $(\tilde{y},\lambda \tilde{y})\in T$ such that $\|y\|_{\mathcal{J}}=\|\hat{y}\|_{\mathcal{K}}$.
Thus, we have that $(\tilde{y},0)\in \lambda I-T$
and $\|y\|_{\mathcal{J}}>0$.
Hence, $N(\lambda I-T)\neq \{0\}$.
Therefore, (4.7) holds.

Next, we show that (4.8) holds.
Suppose that $R(\lambda I-T)=L_W^2(\mathcal{J})$.
For any $\tilde{\hat{g}}\in \hat{L}_W^2(\mathcal{J})$,
set $g:=\hat{g}|_{\mathcal{J}^+}$. Then $\tilde{g}\in L_W^2(\mathcal{J})$.
There exists $\tilde{y}\in D(T)$ such that $(\tilde{y},\tilde{g})\in \lambda I-T$,
i.e., $(\tilde{y},\lambda \tilde{y}-\tilde{g})\in T$.
Take $\tilde{\hat{y}}\in \hat{L}_W^2(\mathcal{J})$ satisfying
$\|\hat{y}\|_{\mathcal{K}}=\|y\|_{\mathcal{J}}$.
Then $(\tilde{\hat{y}},\lambda \tilde{\hat{y}}-\tilde{\hat{g}})\in \hat{T}$,
i.e., $(\tilde{\hat{y}},\tilde{\hat{g}})\in \lambda I-\hat{T}$.
This yields that $R(\lambda I-\hat{T})=\hat{L}_W^2(\mathcal{J})$.
On the other hand, suppose that $R(\lambda I-\hat{T})=\hat{L}_W^2(\mathcal{J})$.
For any $\tilde{f}\in L_W^2(\mathcal{J})$,
take $\tilde{\hat{f}}\in \hat{L}_W^2(\mathcal{J})$ satisfying
$\|\hat{f}\|_{\mathcal{K}}=\|f\|_{\mathcal{J}}$.
Then, there exists $\tilde{\hat{x}}\in D(\hat{T})$ such that
$(\tilde{\hat{x}},\tilde{\hat{f}})\in \lambda I-\hat{T}$, i.e.,
$(\tilde{\hat{x}},\lambda \tilde{\hat{x}}-\tilde{\hat{f}})\in \hat{T}$.
Then, by (4.2), there exists $(\tilde{x},\tilde{h})\in T$ satisfying
$\|x\|_{\mathcal{J}}=\|\hat{x}\|_{\mathcal{K}}$
and $\|h\|_{\mathcal{J}}=\|\lambda \hat{x}-\hat{f}\|_{\mathcal{K}}=\|\lambda x-f\|_{\mathcal{J}}$.
Consequently, $\tilde{h}=\lambda \tilde{x}-\tilde{f}$ and $(\tilde{x},\tilde{f})\in \lambda I-T$.
Hence, $R(\lambda I-T)=L_W^2(\mathcal{J})$. Therefore, (4.8) holds.
This completes the proof.
\medskip

\noindent{\bf Remark 4.1.}
(4.1) and (4.2) are called the zero extensions in the discrete case.
This problem can be easily solved in the continuous case but hard in the discrete case.
For the zero extensions and their properties in analogy with those in Propositions 4.1-4.3 in the continuous case,
please see [3, 7].
\medskip

Note that $H_1$ and $H_{1,r}$ are self-adjoint subspaces
in $(L_W^2(\mathcal{I}))^2$ and $(L_W^2(\mathcal{I}_r))^2$, respectively.
It is difficult to study the convergence of $H_{1,r}$ to $H_1$
in some sense since $L_W^2(\mathcal{I})$ and $L_W^2(\mathcal{I}_r)$ are different spaces.
In order to overcome this problem, we respectively extend
$L_W^2(\mathcal{I}_r)$ and $H_{1,r}$ to be $\hat{L}_W^2(\mathcal{I}_r)$ and $\hat{H}_{1,r}$ by (4.1) and (4.2).
Let $P_r$ be the orthogonal projection from $L_W^2(\mathcal{I})$ to $\hat{L}_W^2(\mathcal{I}_r)$.
Define\vspace{-0.2cm}
$$H_{1,r}':=\hat{H}_{1,r}G(P_r).                                            \vspace{-0.2cm}                  \eqno (4.9)$$

The following result gives the relationship between the spectra of $H_{1,r}', \hat{H}_{1,r}$, and $H_1$,
which is a direct consequence of Propositions 4.2, 4.3, and Lemma 2.3.
\medskip

\noindent\textit{{\bf Lemma 4.1.}
Let $H_1$ be an SSE of $H_0$, and $H_{1,r}$ the induced regular SSE of $H_1$ on $\mathcal{I}_r$.
Then $\hat{H}_{1,r}$ and $H_{1,r}'$ are self-adjoint subspaces in $(\hat{L}_W^2(\mathcal{I}_r))^2$ and $(L_W^2(\mathcal{I}))^2$, respectively,
$D(H_{1,r}')=D(\hat{H}_{1,r})\oplus (\hat{L}_W^2(\mathcal{I}_r))^\bot$,
$\sigma(\hat{H}_{1,r})=\sigma(H_{1,r})$,
and $\sigma(H_{1,r}')=\sigma(\hat{H}_{1,r})\cup \{0\}=\sigma(H_{1,r})\cup \{0\}$.}
\medskip

\section{Spectral approximation in the limit circle case}

In this section, we shall study the regular approximation of spectra of (1.1) in the case that
$\mathscr{L}$ is in l.c.c. at $t=+\infty$. In this case, we shall show that
$\{H_{1,r}\}_{r=1}^\infty$ is not only spectrally inclusive but also spectrally exact for any given $H_1$.
In addition, we obtain explicit approximation relations and give their error estimates.
We always assume that $(\mathbf{A}_1)-(\mathbf{A}_3)$ hold in this section.

For convenience, for any $\tilde{y}\in L_W^2(\mathcal{I})$ and for any $y\in \tilde{y}$, denote\vspace{-0.2cm}
$$\tilde{y}^r:=P_r\tilde{y},\;\;\;\;y_r:=y|_{\mathcal{I}_r^+}.                                      \vspace{-0.2cm}\eqno (5.1)$$
Then, $\tilde{y}^r\in \hat{L}_W^2(\mathcal{I}_r)$, $\tilde{y}_r\in {L}_W^2(\mathcal{I}_r)$, and\vspace{-0.2cm}
$$W(t)R(y^r)(t)=W(t)R(y_r)(t)=W(t)R(y)(t),\qquad t\in \mathcal{I}_r.                               \vspace{-0.2cm}\eqno (5.2) $$

\noindent\textit{{\bf Theorem 5.1.} Assume that $\mathscr{L}$ is in l.c.c. at $t=+\infty$.
Let $H_1$ be any fixed SSE of $H_0$, and $H_{1,r}$ the
induced regular SSE of $H_1$ on $\mathcal{I}_r$, where
$H_1$ and $H_{1,r}$ are determined by (3.5) and (3.7), respectively.
And let $H_{1,r}'$ be defined by (4.9). Then
\begin{itemize}\vspace{-0.2cm}
\item[{\rm (i)}]  $\{H_{1,r}'\}$ is SRC to $H_1$;\vspace{-0.2cm}
\item[{\rm (ii)}] $\{H_{1,r}\}$ is spectrally inclusive for $H_1$ if $0\not\in \sigma(H_1)$.
\end{itemize}}

\noindent{\bf Proof.}
The proof of assertion (i) is divided into three steps:

{\bf Step 1}. Construct a core of $H_1$.

Let\vspace{-0.2cm}
$$C(H_1)=H_{00}\dotplus L\{{\beta}_1,\ldots,{\beta}_{2n}\},    \vspace{-0.1cm}                                     \eqno (5.3) $$
where ${\beta}_i=(\omega_i,\tilde{\tau}_i)\in H, 1\leq i\leq 2n$, are given by (3.4).
By the discussion for Case 1 in Section 3,
$\{{\beta}_1,\ldots,{\beta}_{2n}\}$ is a GKN-set for $\{H_0,H_0^*\}$.
By [29, Theorem 4.2] one gets that\vspace{-0.2cm}
$$H_1=H_0\dotplus L\{{\beta}_1,\ldots,{\beta}_{2n}\},      \vspace{-0.2cm}                        $$
which, together with the fact that $\overline{H}_{00}=H_0$, implies that $C(H_1)$ is a core of $H_1$.

{\bf Step 2}. For any $(y,\tilde{g})\in C(H_1)$, there exists $r_0\in \mathbf{Z^+}$ such that
$(\tilde{y},\tilde{g}^r)\in H_{1,r}'$ for all $r\geq r_0$.

In order to show that this assertion holds, it suffices to show that for any $(y,\tilde{g})\in C(H_1)$,
there exists $r_0\in \mathbf{Z^+}$ such that $(y_r,\tilde{g}_r)\in H_{1,r}$ for all $r\geq r_0$.
In fact, for each $(y,\tilde{g})\in C(H_1)$, if $(y_r,\tilde{g}_r)\in H_{1,r}$, then $(\tilde{y}^r,\tilde{g}^r)\in \hat{H}_{1,r}$.
In addition, since $(\tilde{y},\tilde{y}^r)\in G(P_r)$, we have that $(\tilde{y},\tilde{g}^r)\in H_{1,r}'$
by the definition of $H_{1,r}'$.

Note that $(\mathbf{A}_2)$ for (1.1) on $\mathcal{I}_r$ holds since $b_r> t_0$.
For any given $(y,\tilde{g})\in H_{00}$, by the definition of $H_{00}$,
there exists $r_0\in \mathbf{Z^+}$ such that
$(y_r,\tilde{g}_r)\in H_{1,r}$ for all $r\geq r_0$.
So it is only needed to show that for any
$(y,\tilde{g})\in L\{\beta_1,\ldots,\beta_{2n}\}$,
$(y_r,\tilde{g}_r)\in H_{1,r}$.
For any given $(y,\tilde{g})\in L\{\beta_1,\ldots,\beta_{2n}\}$,
there exist $d_1,\ldots,d_{2n}\in \mathbf{C}$ such that\vspace{-0.2cm}
$$y=\sum_{i=1}^{2n}d_i\omega_i,\;\;\;\;\tilde{g}=\sum_{i=1}^{2n}d_i\tilde{\tau}_i.\vspace{-0.2cm}$$
Since $(y,\tilde{g})\in H_1$, by (3.5) we get that\vspace{-0.2cm}
$$(y,\omega_i)(a)-(y,\omega_i)(+\infty)=0,\;\;1\leq i\leq 2n.                              \vspace{-0.2cm}                \eqno (5.4)$$
In addition, since $\omega_i$, $1\leq i\leq 2n$, are
solutions of $(1.1_\lambda)$ with $\lambda\in \mathbf{R}$ on $[t_0+1,+\infty)$ by (3.4),
it follows from (2.2) that\vspace{-0.2cm}
$$(y,\omega_i)(t)=(y,\omega_i)(+\infty)\;\;{\rm for}\;t\geq t_0+1,\;\;1\leq i\leq 2n.               \vspace{-0.2cm}      \eqno (5.5)$$
Noting that $b_r>t_0$, by (5.4)-(5.5) one has that\vspace{-0.2cm}
$$(y,\omega_i)(a)-(y,\omega_i)(b_r+1)=0,\;\;1\leq i\leq 2n,\vspace{-0.2cm}$$
which yields that $(y_r,\tilde{g}_r)\in H_{1,r}$ by (3.7).
Hence, the assertion in this step holds.

{\bf Step 3}.  $\{H_{1,r}'\}$ and $H_1$ satisfy the conditions in Lemma 2.9.

It follows from Lemma 4.1 that $H_{1,r}', r\geq 1$, are self-adjoint subspaces in $(L_W^2(\mathcal{I}))^2$.
By the assertion in Step 2, we get that
for any $(y,\tilde{g})\in C(H_1)$, there exists a $r_0\in \mathbf{Z^+}$ such that
$(\tilde{y},\tilde{g}^r)\in H_{1,r}'$ for $r\geq r_0$. Since $\|y\|=\|\tilde{y}\|$ and $\tilde{g}=\lim_{r\to\infty} \tilde{g}^r$,
all the conditions in Lemma 2.9 are satisfied. Therefore,
$\{H_{1,r}'\}$ is SRC to $H_1$ by Lemma 2.9.

Assertion (ii) can be directly derived from assertion (i) and Lemmas 2.10 and 4.1. This completes the proof.\medskip

Next, in order to show that $\{H_{1,r}\}$ is spectrally exact for $H_1$,
we shall give the explicit representations of the resolvents
of $H_1$ and $H_{1,r}$ in terms of the Green functions,
respectively, which will play an important role in the discussion of
norm resolvent convergence (for the concept of norm resolvent convergence for self-adjoint subspaces, please see [30, Definition 4.1]),
spectral exactness, and some other topics.
\medskip

\noindent\textit{{\bf Proposition 5.1.} Assume that $\mathscr{L}$ is in l.c.c. at $t=+\infty$.
Let $H_1$ be any SSE of $H_0$.
For any $z\in \rho(H_1)$,
let $\Phi(t,z)=(\phi_1,\ldots,\phi_{2n})(t,z)$ be a standard fundamental solution matrix of $(1.1_z)$ with $\Phi(a,z)=I_{2n}$.
Then, for any $\tilde{g}\in L_W^2(\mathcal{I})$,\vspace{-0.1cm}
$$y(t)=(z I-H_1)^{-1}(\tilde{g})(t)=\sum\limits^{+\infty}_{s=a}G(t,s,z)W(s)R(g)(s),\;\;t\in \mathcal{I},   \vspace{-0.1cm}\eqno (5.6)$$
where\vspace{-0.2cm}
$$G(t,s,z)=\left\{ \begin{array}{cc}
                   \Phi(t,z)M_0R(\Phi)^*(s,\bar{z}), & a\leq s< t< +\infty,\\
                   \Phi(t,z)N_0R(\Phi)^*(s,\bar{z}), & a\leq t\leq s< +\infty,
              \end{array}\!\right.                                                                         \vspace{-0.1cm}\eqno (5.7)$$
is called the Green function of the resolvent $(z I-H_1)^{-1}$, while $M_0$ and $N_0$
are determined by (5.9), (5.11), and (5.12).}
\medskip

\noindent{\bf Proof.} For any fixed $z\in \rho(H_1)$ and for any given $\tilde{g}\in L_W^2(\mathcal{I})$,
from $y=(zI-H_1)^{-1}\tilde{g}$, one has that
$(y,\tilde{g})\in zI-H_1$, and thus $(y,z\tilde{y}-\tilde{g})\in H_1$,
which implies that
$$\mathscr{L}(y)(t)=W(t)R(zy-g)(t),\,\,t\in \mathcal{I};\vspace{-0.2cm}$$
that is,\vspace{-0.2cm}
$$J\Delta y(t)-P(t)R(y)(t)=zW(t)R(y)(t)-W(t)R(g)(t),\,\,t\in \mathcal{I}.   \vspace{-0.2cm}              $$
By the variation of constants formula, every solution $y$ can be given by\vspace{-0.2cm}
$$y(t)=\Phi(t,z)y(a)+\Phi(t,z)J
\sum\limits^{t-1}_{s=a}R(\Phi)^*(s,\bar{z})W(s)R(g)(s),\;t\in \mathcal{I}, \vspace{-0.2cm}                               \eqno (5.8)$$
where we promise that
$$\sum\limits^{a-1}_{s=a}R(\Phi)^*(s,\bar{z})W(s)R(g)(s)=0.$$
Denote
$$\Omega(t):=(\omega_1,\ldots,\omega_{2n})(t),                                                                            \eqno (5.9)$$
where $\omega_i, 1\leq i\leq 2n$, are defined by (3.4).
In view of $y\in D(H_1)$, it follows from (3.5) that\vspace{-0.2cm}
$$\Omega^*(a)Jy(a)=\lim_{t\rightarrow +\infty}\Omega^*(t)Jy(t).     \vspace{-0.2cm}                                   \eqno (5.10)$$
Inserting (5.8) into (5.10), we get that\vspace{-0.2cm}
$$\begin{array}{rrll}
& &\Omega^*(a)Jy(a)\\[0.5ex]
&=&\lim_{t\to +\infty}\Big{\{}\Big{[}\Omega^*(t)J\Phi(t,z)\Big{]}y(a)
+\Big{[}\Omega^*(t)J\Phi(t,z)\Big{]}
J\sum\limits^{t-1}_{s=a}R(\Phi)^*(s,\bar{z})W(s)R(g)(s)\Big{\}}.
\end{array}\vspace{-0.2cm}$$
Since $\mathscr{L}$ is in l.c.c. at $t=+\infty$, we get that
$\phi_1,\ldots,\phi_{2n}\in \mathcal{L}_W^2(\mathcal{I})$, which, together with
(2.1) and $g\in \mathcal{L}_W^2(\mathcal{I})$,
yields that\vspace{-0.2cm}
$$\lim_{t\rightarrow +\infty}
\sum\limits^{t-1}_{s=a}R(\Phi)^*(s,\bar{z})W(s)R(g)(s)
=\sum\limits^{+\infty}_{s=a}R(\Phi)^*(s,\bar{z})W(s)R(g)(s),\vspace{-0.2cm}$$
and\vspace{-0.2cm}
$$K:=\lim_{t\rightarrow +\infty}\Omega^*(t)J\Phi(t,z)   \vspace{-0.2cm}                                            \eqno (5.11)$$
exist and are finite.
So, we have\vspace{-0.2cm}
$$\big(\Omega^*(a)J-K\big)y(a)
=KJ\sum\limits^{+\infty}_{s=a}R(\Phi)^*(s,\bar{z})W(s)R(g)(s).\vspace{-0.2cm}$$
By the fact that $z\in \rho(H_1)$, it can be easily verified that
$\Omega^*(a)J-K$ is invertible.
Thus,\vspace{-0.2cm}
$$y(a)=(\Omega^*(a)J-K)^{-1}KJ\sum\limits^{+\infty}_{s=a}R(\Phi)^*(s,\bar{z})W(s)R(g)(s).\vspace{-0.2cm}$$
Inserting it into (5.8), we get that\vspace{-0.2cm}
$$\begin{array}{rrll}
y(t)&\!=\!&\Phi(t,z)(\Omega^*(a)J-K)^{-1}KJ\sum\limits^{+\infty}_{s=a}R(\Phi)^*(s,\bar{z})W(s)R(g)(s)\\
&\!\!&+\Phi(t,z)J\sum\limits^{t-1}_{s=a}R(\Phi)^*(s,\bar{z})W(s)R(g)(s)\\
&\!=\!&\Phi(t,z)M_0\sum\limits^{t-1}_{s=a}R(\Phi)^*(s,\bar{z})W(s)R(g)(s)
\!+\!\Phi(t,z)N_0\sum\limits^{+\infty}_{s=t}R(\Phi)^*(s,\bar{z})W(s)R(g)(s),
\end{array}\vspace{-0.1cm}$$
where\vspace{-0.2cm}
$$M_0=N_0+J,\;\;\;\;N_0=(\Omega^*(a)J-K)^{-1}KJ.     \vspace{-0.2cm}                                              \eqno (5.12)$$
Therefore, we can write\vspace{-0.2cm}
$$y(t)=\sum\limits^{+\infty}_{s=a}G(t,s,z)W(s)R(g)(s)\vspace{-0.2cm},$$
where $G(t,s,z)$ is specified in (5.7).
This completes the proof.\medskip

Next, consider the explicit representation of resolvent $(z I-H_{1,r})^{-1}$.
With a similar argument in the proof of Proposition 5.1, one can easily show the following result:\medskip

\noindent\textit{{\bf Proposition 5.2.} For any $z\in \rho(H_{1,r})$ and for any $\tilde{g}\in L_W^2(\mathcal{I}_r)$,\vspace{-0.2cm}
$$y(t)=(z I-H_{1,r})^{-1}(\tilde{g})(t)=\sum\limits^{b_r}_{s=a}G_r(t,s,z)W(s)R(g)(s),\;\;t\in \mathcal{I}_r,\vspace{-0.2cm}\eqno (5.13)$$
where\vspace{-0.2cm}
$$G_r(t,s,z)=\left\{ \begin{array}{cc}
                \Phi(t,z)M_rR(\Phi)^*(s,\bar{z}), & a\leq s< t\leq b_r,\\
                \Phi(t,z)N_rR(\Phi)^*(s,\bar{z}), & a\leq t\leq s\leq b_r,
              \end{array}\right.         \vspace{-0.1cm}                                                           \eqno (5.14)$$
is called the Green function of the resolvent $(z I-H_{1,r})^{-1}$, while
$M_r$ and $N_r$ are determined by\vspace{-0.1cm}
$$
M_r=N_r+J,\;\;\;N_r=(\Omega^*(a)J-K_r)^{-1}K_rJ,\;\;\;K_r=\Omega^*(b_r+1)J\Phi(b_r+1,z),\vspace{-0.1cm}           \eqno (5.15)$$}
where $\Omega(t)$ is specified by (5.9), and $\Phi(t,z)$ is specified in Proposition 5.1.
\medskip

Let $A=(a_{ij})\in \mathbf{C}^{k\times l}$ and $\xi=(\xi_1,\ldots,\xi_l)^{\rm T}\in \mathbf{C}^{l}$.
Define their norms as\vspace{-0.1cm}
$$\|A\|_1:=(\sum\limits^{k}_{i=1}\sum\limits^{l}_{j=1}|a_{ij}|^2)^{1/2},\;\;\;\;
\|\xi\|_1=(\sum\limits^{l}_{j=1}|\xi_j|^2)^{1/2}.\vspace{-0.1cm}$$
Then\vspace{-0.1cm}
$$\|A\xi\|_1\leq \|A\|_1\|\xi\|_1,\;\;\;\;\|AB\|_1\leq \|A\|_1\|B\|_1,\;\;\forall \; B\in \mathbf{C}^{l\times n}.\vspace{-0.1cm}$$

It follows from (5.11) and (5.15) that $K_r\rightarrow K$ as $r\rightarrow +\infty$.
So one can get the following result by Propositions 5.1 and 5.2.
\medskip

\noindent\textit{{\bf Proposition 5.3.} $M_r\rightarrow M_0$ and $N_r\rightarrow N_0$ as $r\rightarrow +\infty$.}
\medskip

Now, we can give the following result about spectral exactness.\medskip

\noindent\textit{{\bf Theorem 5.2.} Assume that $\mathscr{L}$ is in l.c.c. at $t=+\infty$.
Let $H_1$ be any fixed SSE of $H_0$, and $H_{1,r}$ the
induced regular SSE of $H_1$ on $\mathcal{I}_r$, where
$H_1$ and $H_{1,r}$ are determined by (3.5) and (3.7), respectively.
And let $\hat{H}_{1,r}$ be defined by (4.2). Then
\begin{itemize}\vspace{-0.2cm}
\item[{\rm (i)}] for any $z\in \rho(H_1)\cap \rho(\hat{H}_{1,r})$,
$\{(zI-\hat{H}_{1,r})^{-1}G(P_r)\}\stackrel{n}{\to} (zI-H_1)^{-1}$;\vspace{-0.2cm}
\item[{\rm (ii)}] $\{H_{1,r}\}$ is spectrally exact for $H_1$ if $0\not\in \si(H_1)$.
\end{itemize}}

\noindent{\bf Proof.}
We first show that assertion (i) holds.
Let $z\in \rho(H_1)\cap \rho(\hat{H}_{1,r})$.
Note that $(zI-\hat{H}_{1,r})^{-1}$ is an operator. So we write $\{(zI-\hat{H}_{1,r})^{-1}G(P_r)\}$
as $\{(zI-\hat{H}_{1,r})^{-1}P_r\}$ for short.
It follows from (5.1)-(5.2) that for any given $\tilde{g}\in L_W^2(\mathcal{I})$,\vspace{-0.1cm}
$$\begin{array}{rrll}
W(t)R((z I-\hat{H}_{1,r})^{-1}P_r\tilde{g})(t)&=&W(t)R((z I-\hat{H}_{1,r})^{-1}\tilde{g}^r)(t)\\
&=&W(t)R((z I-H_{1,r})^{-1}\tilde{g}_r)(t),\;\;t\in \mathcal{I}_r.
\end{array}\vspace{-0.1cm}$$
It yields that\vspace{-0.2cm}
$$\delta_r(\tilde{g}):=\|(z I-H_1)^{-1}\tilde{g}-(z I-\hat{H}_{1,r})^{-1}P_r\tilde{g}\|^2
=\delta_{r1}(\tilde{g})+\delta_{r2}(\tilde{g}),\vspace{-0.2cm}                                                     \eqno (5.16)$$
where\vspace{-0.2cm}
$$\delta_{r1}(\tilde{g})\!\!:=\!\!\sum\limits^{b_r}_{t=a}R((z I\!-\!H_1)^{-1}\tilde{g}\!-\!(z I\!-\!H_{1,r})^{-1}\tilde{g}_r)^*(t)W(t)
R((z I\!-\!H_1)^{-1}\tilde{g}\!-\!(z I\!-\!H_{1,r})^{-1}\tilde{g}_r)(t),\vspace{-0.2cm}$$
\vspace{-0.2cm}
$$\delta_{r2}(\tilde{g})\!\!:=\!\!\!\sum\limits^{+\infty}_{t=b_r+1}R((z I-H_1)^{-1}\tilde{g})^*(t)W(t)R((z I-H_1)^{-1}\tilde{g})(t).\vspace{-0.1cm}$$

Now, we first consider $\delta_{r1}(\tilde{g})$.
By Propositions 5.1 and 5.2 we get that\vspace{-0.1cm}
$$\begin{array}{rrll}
&&(z I-H_1)^{-1}\tilde{g}(t)-(z I-H_{1,r})^{-1}\tilde{g}_r(t)\\
&=&\sum\limits^{+\infty}_{s=a}G(t,s,z)W(s)R(g)(s)-\sum\limits^{b_r}_{s=a}G_r(t,s,z)W(s)R(g_r)(s)\\
&=&T_{r1}(\tilde{g})(t)+T_{r2}(\tilde{g})(t)+T_{r3}(\tilde{g})(t),\;t\in \mathcal{I}_r,
\end{array}\vspace{-0.1cm}$$
where\vspace{-0.1cm}
$$\begin{array}{rrll}
T_{r1}(\tilde{g})(t):&=&\Phi(t,z)(M_0-M_r)\sum\limits^{t-1}_{s=a}R(\Phi)^*(s,\bar{z})W(s)R(g)(s),\\[1.2ex]                                                                                                                  T_{r2}(\tilde{g})(t):&=&\Phi(t,z)(N_0-N_r)\sum\limits^{b_r}_{s=t}R(\Phi)^*(s,\bar{z})W(s)R(g)(s),\\[1.2ex]                                                                                                                  T_{r3}(\tilde{g})(t):&=&\Phi(t,z)N_0\sum\limits^{+\infty}_{s=b_r+1}R(\Phi)^*(s,\bar{z})W(s)R(g)(s).                                                                                                                   \end{array}\vspace{-0.2cm}$$
Consequently,\vspace{-0.2cm}
$$\delta_{r1}(\tilde{g})\leq 3\sum^{3}_{i=1}\sum_{t\in \mathcal{I}_r}R(T_{ri}(\tilde{g}))^*(t)W(t)R(T_{ri}(\tilde{g}))(t).\vspace{-0.2cm}\eqno (5.17)$$
Denote\vspace{-0.2cm}
$$\begin{array}{cc}
    m_0:=\|M_0\|_1, \;\;\;n_0:=\|N_0\|_1, \;\;\; m_r:=\|M_0-M_r\|_1,\;\;\;n_r:=\|N_0-N_r\|_1,\\[0.8ex]
\alpha_0(z):=\max\limits_{1\leq i\leq 2n}\{\|\phi_i(\cdot,z)\|\}, \;
\alpha_r(z):=\max\limits_{1\leq i\leq 2n}\{\sum\limits^{\infty}_{t=b_r+1}R(\phi_i)^*(t,z)W(t)R(\phi_i)(t,z)\}.
\end{array}\vspace{-0.2cm}                                                                                                    \eqno (5.18)$$
Then $m_r\rightarrow 0$ and $n_r\rightarrow 0$ by Proposition 5.3, and $\alpha_r(z)\rightarrow 0$ as $r\rightarrow \infty$.
For convenience, denote\vspace{-0.2cm}
$$h_r(t):=(M_0-M_r)\sum\limits^{t-1}_{s=a}R(\Phi)^*(s,\bar{z})W(s)R(g)(s),\vspace{-0.2cm}$$
then\vspace{-0.2cm}
$$R(T_{r1}(\tilde{g}))(t)=R(\Phi(t,z)h_r(t))=R(\Phi(t,z))h_r(t)+{\rm diag}\{I_n,0\}\Phi(t+1,z)\Delta h_r(t),     \vspace{-0.2cm} \eqno (5.19)$$
$$\|h_r(t)\|_1\leq \|M_0-M_r\|_1\|\sum\limits^{t-1}_{s=a}R(\Phi)^*(s,\bar{z})W(s)R(g)(s)\|_1\leq \sqrt{2n}m_r\alpha_0(\bar{z})\|\tilde{g}\|.
                                                                                                            \vspace{-0.1cm}\eqno (5.20)$$
In addition, it follows from $(1.1_z)$ that\vspace{-0.1cm}
$$\Phi(t+1,z)=\left(\begin{array}{cc}
                      I_n & 0 \\
                      C(t)-zW_1(t) & I_n-A^*(t)
                    \end{array}\right)R(\Phi)(t,z).\vspace{-0.1cm}$$
Inserting it into (5.19), we get that\vspace{-0.2cm}
$$R(T_{r1}(\tilde{g}))(t)=R(\Phi(t,z))h_r(t)+{\rm diag}\{I_n,0\}R(\Phi)(t,z)\Delta h_r(t).\vspace{-0.1cm}$$
Therefore,\vspace{-0.1cm}
$$\begin{array}{rrll}
&&R(T_{r1}(\tilde{g}))^*(t)W(t)R(T_{r1}(\tilde{g}))(t)\\[1.2ex]
&=&h_r^*(t)R(\Phi)^*(t,z)W(t)R(\Phi)(t,z)h_r(t)\\[1.2ex]
&&+h_r^*(t)R(\Phi)^*(t,z){\rm diag}\{W_1(t),0\}R(\Phi)(t,z)\Delta h_r(t)\\[1.2ex]
&&+\Delta h_r^*(t)R(\Phi)^*(t,z){\rm diag}\{W_1(t),0\}R(\Phi)(t,z)h_r(t)\\[1.2ex]
&&+\Delta h_r^*(t)R(\Phi)^*(t,z){\rm diag}\{W_1(t),0\}R(\Phi)(t,z)\Delta h_r(t).
\end{array} \vspace{-0.2cm}                                                                                         \eqno (5.21)$$
Since $\mathscr{L}$ is in l.c.c. at $t=+\infty$, all the solutions of $(1.1_z)$ are in $\mathcal{L}_W^2(\mathcal{I})$,
and so $\Phi(\cdot,z)\in \mathcal{L}_W^2(\mathcal{I})$. It follows that all the diagonal entries of $R(\Phi)^*(t,z)W(t)R(\Phi)(t,z)$
are nonnegative and absolutely summable over $[a,+\infty)$.
In addition, using the nonnegativity of $W(t)$, one has that\vspace{-0.1cm}
$$\sum\limits^{+\infty}_{t=a}\|R(\Phi)^*(t,z)W(t)R(\Phi)(t,z)\|_1 \leq  2n{\alpha}_0^2(z),\vspace{-0.1cm}     \eqno (5.22)$$
which implies that\vspace{-0.1cm}
$$\sum\limits^{+\infty}_{t=a}\|R(\Phi)^*(t,z){\rm diag}\{W_1(t),0\}R(\Phi)(t,z)\|_1 \leq 2n{\alpha}_0^2(z),\vspace{-0.1cm}$$
which, together with (5.20)-(5.22), yields that\vspace{-0.1cm}
$$\begin{array}{rrll}
&&\sum\limits_{t\in \mathcal{I}_r}R(T_{r1}(\tilde{g}))^*(t)W(t)R(T_{r1}(\tilde{g}))(t)\\
&\leq& \sum\limits_{t\in \mathcal{I}_r}\|R(\Phi)^*(t,z)W(t)R(\Phi)(t,z)\|_1
[\|h_r(t)\|_1^2+2\|h_r(t)\|_1\|\Delta h_r(t)\|_1+\|\Delta h_r(t)\|_1^2]\\
&\leq& 36n^2{\alpha}_0^2(z){\alpha}_0^2(\bar{z}){m}_r^2\|\tilde{g}\|^2.
\end{array}                                                                                         \vspace{-0.1cm}\eqno (5.23)$$
Similarly, we get that\vspace{-0.1cm}
$$\begin{array}{rrll}
&&\sum\limits_{t\in \mathcal{I}_r}R(T_{r2}(\tilde{g}))^*(t)W(t)R(T_{r2}(\tilde{g}))(t)\leq
36n^2{\alpha}_0^2(z){\alpha}_0^2(\bar{z}){n}_r^2\|\tilde{g}\|^2,\\[2.6ex]
&&\sum\limits_{t\in \mathcal{I}_r}R(T_{r3}(\tilde{g}))^*(t)W(t)R(T_{r3}(\tilde{g}))(t)\leq
4n^2{\alpha}_0^2(z){\alpha}_r(\bar{z}){n}_0^2\|\tilde{g}\|^2.
\end{array}\vspace{-0.1cm}$$
Thus, from (5.17), it follows that\vspace{-0.1cm}
$$\delta_{r1}(\tilde{g})\leq 12n^2{\alpha}_0^2(z)
\Big{[}9{\alpha}_0^2(\bar{z})({m}_r^2+{n}_r^2)+{\alpha}_r(\bar{z}){n}_0^2\Big{]}\|\tilde{g}\|^2.     \vspace{-0.1cm}\eqno(5.24)$$

With a similar argument to that used for $\delta_{r1}(\tilde{g})$,
one can show that\vspace{-0.2cm}
$$\delta_{r2}(\tilde{g})\leq 72n^2({m}_0^2+{n}_0^2){\alpha}_0^2(\bar{z}){\alpha}_r(z)\|\tilde{g}\|^2,\vspace{-0.2cm} $$
which, together with (5.24), yields that\vspace{-0.2cm}
$$\delta_{r}(\tilde{g})=\delta_{r1}(\tilde{g})+\delta_{r2}(\tilde{g})\leq \eta(r)\|\tilde{g}\|^2,\vspace{-0.2cm}   \eqno (5.25)$$
where\vspace{-0.2cm}
$$\eta(r)=12n^2\Big{[}9\alpha_0^2(z){\alpha}_0^2(\bar{z})({m}_r^2+{n}_r^2)+
{n}_0^2\alpha_0^2(z){\alpha}_r(\bar{z})+
6({m}_0^2+{n}_0^2)\alpha_0^2(\bar{z}){\alpha}_r(z)\Big{]}.\vspace{-0.2cm}                                         \eqno (5.26)$$
It is obvious that $\eta(r)\rightarrow 0$ as $r\rightarrow \infty$,
which implies that assertion (i) holds.

It can be directly derived from Lemma 4.1, Theorem 5.1, assertion (i), and Lemma 2.10 that
$\{{ H}_{1,r}\}$ is spectrally exact for $H_1$ by the assumption that $0\not\in \si(H_1)$.
The whole proof is complete.\medskip

In order to further study how to approximate the spectrum of $H_1$ by those of $\{H_{1,r}\}$, we first give the following
useful result:\medskip

\noindent\textit{{\bf Theorem 5.3.}
Every self-adjoint subspace extension $H_1$ of $H_0$ has a pure discrete spectrum in the case that
$\mathscr{L}$ is in l.c.c. at $t=+\infty$.}
\medskip

\noindent{\bf Proof.} According to [42, Theorems 6.7 and 6.10] and Lemma 2.1, it suffices to prove
that $(zI-H_1)^{-1}$ is a Hilbert-Schmidt operator for any $z\in \rho(H_1)$.

By Proposition 5.1, for any $z\in \rho(H_1)$ and any $\tilde{g}\in L_W^2(\mathcal{I})$,\vspace{-0.2cm}
$$(z I-H_1)^{-1}(\tilde{g})(t)=\sum\limits^{+\infty}_{s=a}G(t,s,z)W(s)R(g)(s),\,\,t\in \mathcal{I},\vspace{-0.1cm}$$
where $G(t,s,z)$ is given by (5.7). Define\vspace{-0.1cm}
$$\begin{array}{rrll}
\mathcal{F}_1(\tilde{g})(t)&:=&\sum\limits^{+\infty}_{s=a}F_1(t,s,z)W(s)R(g)(s),\,\,t\in \mathcal{I},\\[0.5ex]
\mathcal{F}_2(\tilde{g})(t)&:=&\sum\limits^{+\infty}_{s=a}F_2(t,s,z)W(s)R(g)(s),\,\,t\in \mathcal{I},
\end{array}        \vspace{-0.1cm} $$
where
$$F_1(t,s,z)=\left\{ \begin{array}{cc}
                   \Phi(t,z)M_0R(\Phi)^*(s,\bar{z}), & a\leq s< t< +\infty,\\
                   0, & a\leq t\leq s< +\infty,
              \end{array}\!\right.                                                                                    $$
$$F_2(t,s,z)=\left\{ \begin{array}{cc}
                   0, & a\leq s< t< +\infty,\\
                   \Phi(t,z)N_0R(\Phi)^*(s,\bar{z}), & a\leq t\leq s< +\infty.
              \end{array}\!\right.                                                                                   $$
Obviously, $(z I-H_1)^{-1}=\mathcal{F}_1+\mathcal{F}_2$. Therefore, it is sufficient to prove that
$\mathcal{F}_1$ and $\mathcal{F}_2$ are both Hilbert-Schmidt operators.
Denote $N:={\rm dim}L_W^2(\mathcal{I})$. In the case of $N<\infty$,
$\mathcal{F}_1$ and $\mathcal{F}_2$ are obviously Hilbert-Schmidt operators.
So, it is only needed to show that this assertion holds in the case of $N=\infty$.
We first prove that $\mathcal{F}_1$ is a Hilbert-Schmidt operator in this case.
Let $\{\tilde{e}_j:j\in \mathbf{Z}^{+}\}$ be an orthonormal basis of $L_W^2(\mathcal{I})$.
Then\vspace{-0.1cm}
$$
\mathcal{F}_1(\tilde{e}_j)(t)=\sum\limits^{+\infty}_{s=a}F_1(t,s,z)W(s)R(e_j)(s)=\Phi(t,z)u_j(t),\,\,t\in \mathcal{I},
\vspace{-0.1cm}$$
where\vspace{-0.1cm}
$$u_j(t):=M_0\sum\limits^{t-1}_{s=a}R(\Phi)^*(s,\bar{z})W(s)R(e_j)(s).\vspace{-0.1cm}$$
Then\vspace{-0.1cm}
$$\|u_j(t)\|_1\leq \|M_0\|_1\|\sum\limits^{t-1}_{s=a}R(\Phi)^*(s,\bar{z})W(s)R(e_j)(s)\|_1
\leq m_0\Big{(}\sum\limits^{2n}_{i=1}|\langle\phi_i(\cdot,\bar{z}),e_j\rangle|^2\Big{)}^{1/2}.
\vspace{-0.1cm}                                                                                                              \eqno (5.27)$$
Similar to the discussions for (5.21) and (5.23) with replacing $h_r(t)$ by $u_j(t)$, one has that
$$\begin{array}{rrll}
\sum\limits^{\infty}_{j=1}\|\mathcal{F}_1(\tilde{e}_j)\|^2
&=&\sum\limits^{\infty}_{j=1}\sum\limits^{+\infty}_{t=a}R(\Phi(t,z)u_j(t))^*W(t)R(\Phi(t,z)u_j(t))\\[2.6ex]
&\leq& 18nm_0^2\alpha_0^2(z)\sum\limits^{2n}_{i=1}\sum\limits^{\infty}_{j=1}|\langle\phi_i(\cdot,\bar{z}),e_j\rangle|^2\\[2.6ex]
&\leq& 36n^2m_0^2\alpha_0^2(z)\alpha_0^2(\bar{z})<\infty,
\end{array}$$
in which (5.27), (5.22), and Parseval's identity have been used.
Therefore, $\mathcal{F}_1$ is a Hilbert-Schmidt operator. Similarly, one can show that
$\mathcal{F}_2$ is a Hilbert-Schmidt operator and thus $(z I-H_1)^{-1}$ is a Hilbert-Schmidt operator.
The proof is complete.\medskip

\noindent{\bf Remark 5.1.}
With a similar argument, by applying the Green function of $H_{1,r}$ given in Proposition 5.2,
it can be easily verified that the resolvent of $H_{1,r}$ is a Hilbert-Schmidt operator.
Hence, the resolvent of $\hat{H}_{1,r}$ is also a Hilbert-Schmidt operator by (4.2).
In addition, denote $\sigma(H_1)\backslash \{0\}:=\{\lambda_{k}:\,k\in \Upsilon\}$,
where $\Upsilon$ denotes the eigenvalue index set of $H_1$.
Then we can further get that
$\sum\limits_{k\in \Upsilon}|\lambda_k|^{-2}<\infty$
by the fact that the resolvent of $H_{1}$ is a Hilbert-Schmidt operator.
\medskip

By Theorem 5.3, $H_1$ has a discrete spectrum in the case that $\mathscr{L}$ is in l.c.c. at $t=+\infty$.
By translating it if necessary, we may suppose that 0 is not an eigenvalue of $H_1$.
The eigenvalues of $H_1$ may be ordered as (multiplicity included):
$$\cdots\leq \lambda_{-2} \leq \lambda_{-1}<0<\lambda_{1} \leq \lambda_{2}\leq\cdots\leq \lambda_{k}\leq\cdots.$$
For convenience, by $\Upsilon$ denote the eigenvalue index set of $H_1$ and $\sigma(H_1)=\{\lambda_{k}:\,k\in \Upsilon\}$.
By (ii) of Theorem 5.2, $\{H_{1,r}\}$ is spectrally exact for $H_1$ if $0\notin \sigma(H_1)$.
Hence, since $0\notin \sigma(H_1)$, there exists $r_0$ such that $0\notin \sigma(H_{1,r})$ for all $r\geq r_0$.
Therefore, for $r\geq r_0$, the eigenvalues of $H_{1,r}$ may be ordered as (multiplicity included):
$$\lambda_{-m(r)}^{(r)} \leq\cdots\leq \lambda_{-2}^{(r)} \leq \lambda_{-1}^{(r)}<0<\lambda_{1}^{(r)}
\leq \lambda_{2}^{(r)}\leq\cdots \leq \lambda_{n(r)}^{(r)},                                       $$
where $m(r)$ and $n(r)$ are the numbers of negative and positive eigenvalues of $H_{1,r}$, respectively.
For convenience, we briefly denote the eigenvalue index set of $H_{1,r}$ by $\Upsilon_r$, and then
$\sigma(H_{1,r})=\{\lambda_{k}^{(r)}:\,k\in \Upsilon_r\}$.
By Lemma 4.1, $\sigma(H_{1,r})=\sigma(\hat{H}_{1,r})$, which implies that $0\in \rho(\hat{H}_{1,r})$ as $r\geq r_0$.
\medskip

\noindent\textit{{\bf Theorem 5.4.} Assume that $\mathscr{L}$ is in l.c.c. at $t=+\infty$.
For each $k\in \Upsilon$, there exists an $r_k\geq r_0$ such that
$k\in \Upsilon_r$ for $r\geq r_k$, and $\lambda_{k}^{(r)}\to \lambda_{k}$ as $r \to \infty$.}\medskip

\noindent{\bf Proof.}
Let\vspace{-0.1cm}
$$S=(-H_1)^{-1},\;\;S_r=(-\hat{H}_{1,r})^{-1},\;r\geq r_0.\vspace{-0.1cm}$$
Then, according to Lemmas 2.3 and 4.1, the proof of Theorem 5.3, and Remark 5.1,
it follows that $S_rP_r$ and $S$ are both self-adjoint and Hilbert-Schmidt operators in $L_W^2(\mathcal{I})$.
And $\mu_k=-1/{\lambda_{k}}$ for $k\in \Upsilon$ and
$\mu_{k}^{(r)}=-1/{\lambda_{k}^{(r)}}$ for $k\in \Upsilon_r$ are eigenvalues of $S$ and $S_rP_r$, respectively, by Lemmas 2.1 and 2.3.
($S_rP_r$ also has 0 as an eigenvalue of infinite multiplicity. But it is not related to $H_{1,r}$ or $H_1$, and so can be ignored.)
Further, $S_rP_r {\to} S$ in norm as $r \to \infty$ by (i) of Theorem 5.2.
It follows that $S_rP_r {\to} S$ in the norm resolvent sense as $r \to \infty$ according to the proof of [23, Theorem 8.18]
(for the concept of convergence of self-adjoint operators in the norm resolvent sense, please see [23, 42]).
Let $E(S_rP_r,\lambda)$ and $E(S,\lambda)$ be spectral families of $S_rP_r$ and $S$, respectively.
Then, by (b) of [23, Theorem 8.23] one has that
for any $\alpha,\beta \in \mathbf{R}\cap \rho(S)$ with $\alpha<\beta$,\vspace{-0.1cm}
$$\|[E(S_rP_r,\beta)-E(S_rP_r,\alpha)]-[E(S,\beta)-E(S,\alpha)]\| \to 0 \;{\rm as}\;r \to \infty,\vspace{-0.1cm}$$
which, together with [42, Theorem 4.35], yields that\vspace{-0.1cm}
$$\dim R[E(S_rP_r,\beta)-E(S_rP_r,\alpha)]=\dim R[E(S,\beta)-E(S,\alpha)]\vspace{-0.1cm}$$
for all sufficiently large $r$. Hence, for each $k\in \Upsilon$, there exists an $r_k\geq r_0$ such that
$k\in \Upsilon_r$ for $r\geq r_k$.

Next, we show that $\lambda_{k}^{(r)}\to \lambda_{k}$ as $r \to \infty$.
To do so, it suffices to prove that $\mu_{k}^{(r)}\to \mu_{k}$ as $r \to \infty$.
The eigenvalues are described by the Courant-Fischer min-max theorem in the case of ${\rm dim}L_W^2(\mathcal{I})<\infty$
and by a min-max principle according to [26, Section 12.1] in the case of ${\rm dim}L_W^2(\mathcal{I})=\infty$, respectively;
that is,\vspace{-0.1cm}
$$\mu_k=\left\{\begin{array}{cc}
 \min\limits_{V_k}\max\limits_{{x\in {V_k},\atop \|x\|=1}}\langle Sx,x\rangle, & k\in \Upsilon \,{\rm with}\, k>0,\\
 \max\limits_{V_k}\min\limits_{{x\in {V_k},\atop \|x\|=1}}\langle Sx,x\rangle, & k\in \Upsilon \,{\rm with}\, k<0,\\
\end{array}\right.                                                                     \vspace{-0.1cm} \eqno (5.28)$$
where $V_k$ runs through all the $|k|$-dimensional subspaces of $L^2_W(\mathcal{I})$.
For $r\geq r_k$, $\mu_{k}^{(r)}$ is similarly expressed in terms of $\langle S_rP_rx,x\rangle$; that is,\vspace{-0.1cm}
$$\mu_k^{(r)}=\left\{\begin{array}{cc}
 \min\limits_{V_k}\max\limits_{{x\in {V_k},\atop \|x\|=1}}\langle S_rP_rx,x\rangle, & k\in \Upsilon \,{\rm with}\, k>0,\\
 \max\limits_{V_k}\min\limits_{{x\in {V_k},\atop \|x\|=1}}\langle S_rP_rx,x\rangle, & k\in \Upsilon \,{\rm with}\, k<0.\\
\end{array}\right.                                                                           \vspace{-0.1cm} \eqno (5.29)$$

We first consider the case that $k\in \Upsilon$ with $k>0$. Let $r\geq r_k$.
It follows from (5.28)-(5.29) that there exist two
$k$-dimensional subspaces $V_k$ and $\tilde{V}_k$ of $L^2_W(\mathcal{I})$ such that\vspace{-0.1cm}
$$\mu_{k}=\max\limits_{{x\in {V_k},\atop \|x\|=1}}\langle Sx,x\rangle,\;\;\;\;
\mu_{k}^{(r)}=\max\limits_{{x\in {\tilde{V}_k},\atop \|x\|=1}}\langle S_rP_rx,x\rangle.   \vspace{-0.1cm} \eqno (5.30)$$
In addition, there exist $x_1\in \tilde{V}_k$ with $\|x_1\|=1$ and $x_2\in V_k$ with $\|x_2\|=1$ such that\vspace{-0.1cm}
$$\max\limits_{{x\in {\tilde{V}_k},\atop \|x\|=1}}\langle Sx,x\rangle=\langle Sx_1,x_1\rangle,\;\;\;\;
\max\limits_{{x\in {V}_k,\atop \|x\|=1}}\langle S_rP_rx,x\rangle=\langle S_rP_rx_2,x_2\rangle. \vspace{-0.1cm}\eqno (5.31)$$
From (5.28)-(5.31), we have\vspace{-0.1cm}
$$\mu_{k}-\mu_{k}^{(r)}
\leq \max\limits_{{x\in {\tilde{V}_k},\atop \|x\|=1}}\langle Sx,x\rangle-
\max\limits_{{x\in {\tilde{V}_k},\atop \|x\|=1}}\langle S_rP_rx,x\rangle
\leq \langle (S-S_rP_r)x_1,x_1\rangle,\vspace{-0.1cm}$$
$$\mu_{k}-\mu_{k}^{(r)}\geq
\max\limits_{{x\in {V_k},\atop \|x\|=1}}\langle Sx,x\rangle-\max\limits_{{x\in {V_k},\atop \|x\|=1}}\langle S_rP_rx,x\rangle
\geq \langle (S-S_rP_r)x_2,x_2\rangle.\vspace{-0.1cm}$$
Therefore, it follows that\vspace{-0.1cm}
$$\begin{array}{rrll}
&&|\mu_{k}-\mu_{k}^{(r)}|\leq \max\{|\langle (S-S_rP_r)x_1,x_1\rangle|, |\langle (S-S_rP_r)x_2,x_2\rangle|\}\\
&\leq& \|S-S_rP_r\| \to 0\,\,{\rm as}\,r \to \infty.
\end{array}                                                                               \vspace{-0.1cm}\eqno (5.32)$$
Thus, $\mu_{k}^{(r)}\to \mu_{k}$ as $r \to \infty$ for $k\in \Upsilon$ with $k>0$.

Similarly, one can get that $\mu_{k}^{(r)}\to \mu_{k}$ as $r \to \infty$ for $k\in \Upsilon$ with $k<0$.
This completes the proof.
\medskip

At the end of this section, we shall try to give an error estimate for the approximation
of $\lambda_{k}$ by $\lambda_{k}^{(r)}$ for each $k\in \Upsilon$. Obviously, it is very important
in numerical analysis and applications. In order to give error estimates of $\lambda_{k}^{(r)}$ to $\lambda_{k}$, in view of
$\lambda_k=-1/{\mu_{k}}$ and $\lambda_{k}^{(r)}=-1/{\mu_{k}^{(r)}}$, we shall first investigate the
error estimates of $\mu_{k}^{(r)}$ to $\mu_{k}$ for $k\in \Upsilon$ instead.\medskip

In view of the arbitrariness of $\lambda\in \mathbf{R}$ in (2.4),
we might as well take $\lambda=0$ in (2.4) in the rest of this section.
\medskip

\noindent\textit{{\bf Proposition 5.4.} Assume that $\mathscr{L}$ is in l.c.c. at $t=+\infty$.
Then, for each $k\in\Upsilon$ and $r\geq r_k$, where $r_k$ is specified in Theorem $5.4$,\vspace{-0.2cm}
$$|\mu_{k}^{(r)}-\mu_{k}|\leq 2\sqrt{3}n\alpha_0[(6m_0^2+7n_0^2)(\|E(a)\|_1^2+\|E(a)B(a)\|_1^2+n)]^{1/2}\varepsilon_r^{1/2},
                                                                                                            \vspace{-0.2cm} \eqno (5.33)$$
where $\alpha_0:=\alpha_0(0)$, $m_0$, and $n_0$ are constants and given by $(5.18)$,
$E(a)=(I_n-A(a))^{-1}$, and $\varepsilon_r$ is completely determined by
the coefficients of $(1.1)$, more precisely, it is determined by $(5.36)$, $(5.38)$, $(5.40)$, and $(5.42)$.
In addition, $\varepsilon_r\to 0$ as $r\to \infty$.}\medskip

\noindent{\bf Proof.}
In view of $0\in \rho(H_1)\cap \rho(\hat{H}_{1,r})$ as $r\geq r_0$,
it follows from (5.26) with $z=0$ and (5.32) that\vspace{-0.1cm}
$$|\mu_{k}^{(r)}-\mu_{k}|\leq 2\sqrt{3}n\alpha_0[9{\alpha}_0^2({m}_r^2+{n}_r^2)+
(6{m}_0^2+7{n}_0^2){\alpha}_r]^{1/2},\;r\geq r_0,                                                     \vspace{-0.1cm}\eqno (5.34)$$
where $\alpha_r:=\alpha_r(0)$ and $\alpha_0:=\alpha_0(0)$,
and $m_r, n_r, m_0, n_0, \alpha_r(z), \alpha_0(z)$ are specified in (5.18).

By the arbitrariness of $\lambda\in \mathbf{R}$ in (2.4), we take $\lambda=0$ in it.
So it follows from (3.4) that $\omega_1,\ldots,\omega_{2n}$ are solutions of $(1.1_\lambda)$ with $\lambda=0$ in $[t_0+1,+\infty)$,
where $t_0$ is specified by $(\mathbf{A}_2)$.
In addition, since $\phi_1(\cdot,0),\ldots,\phi_{2n}(\cdot,0)$ are solutions of $(1.1_z)$ with $z=0$ in $\mathcal{I}$
and $b_r> t_0$, by (2.2), (5.11), and (5.15)
we get that $K=K_r$, which, together with (5.12), yields that $M_r=M_0$, $N_r=N_0$, and thus $m_r=n_r=0$ by (5.18).

Now, it remains to estimate $\alpha_r$. By $(\mathbf{A}_1)$, we get that every solution $y$ of
$(1.1_z)$ with $z=0$ satisfies\vspace{-0.1cm}
$$ R(y)(t+1)=U(t)R(y)(t),\;\;\;t\geq a,                                                \vspace{-0.1cm}\eqno (5.35)$$
where $E(t)=(I_n-A(t))^{-1}$,\vspace{-0.1cm}
$$U(t)=\left(\begin{array}{cc}
  E(t+1)+E(t+1)B(t+1)C(t) & E(t+1)B(t+1)(I_n-A^*(t)) \\
  C(t) &  I_n-A^*(t)
\end{array}\right).                                                                     \vspace{-0.1cm}\eqno (5.36)$$
It follows that\vspace{-0.1cm}
$$\begin{array}{rrll}
&&R(y)^*(t+1)W(t+1)R(y)(t+1)\\
&=&R(y)^*(t)U^*(t)W(t+1)U(t)R(y)(t)\\
&=&R(y)^*(t-1)U^*(t-1)U^*(t)W(t+1)U(t)U(t-1)R(y)(t-1)\\
&=&R(y)^*(a)V^*(t)W(t+1)V(t)R(y)(a),\;\;t\geq a,
\end{array}                                                                            \vspace{-0.1cm}\eqno (5.37) $$
where\vspace{-0.1cm}
$$V(t):= U(t)U(t-1)\cdots U(a),\;\;t\geq a.                                             \vspace{-0.1cm}\eqno (5.38)$$
Since $\mathscr{L}$ is in l.c.c. at $t=+\infty$, it follows from (5.37) that\vspace{-0.1cm}
$$\begin{array}{rrll}
&&\sum\limits^{\infty}_{t=b_r}R(y)^*(t+1)W(t+1)R(y)(t+1)\\[2.0ex]
&=&R(y)^*(a)D_rR(y)(a)\to 0\;{\rm as}\; r\to \infty,
\end{array}                                                                            \vspace{-0.1cm}\eqno (5.39)$$
where\vspace{-0.1cm}
$$D_r:=\sum\limits^{\infty}_{t=b_r}V^*(t)W(t+1)V(t).                                    \vspace{-0.1cm}\eqno (5.40)$$
Then, $D_r$ is positive semi-definite since $W(t)$ is positive semi-definite for $t\geq a$.
In addition, since\vspace{-0.1cm}
$$R(y)(a)=\left(\begin{array}{cc}
  E(a) & E(a)B(a) \\
  0 &  I_n
\end{array}\right)y(a),                                                                   \vspace{-0.1cm}\eqno (5.41)$$
$R(y)(a)$ can be taken any complex vector belonging to $\mathbf{C}^{2n}$.
Denote\vspace{-0.1cm}
$$\varepsilon_r:=\|D_r\|_1.                                                              \vspace{-0.1cm}\eqno (5.42)$$
Then, combining the positive semi-definiteness of $D_r$ and the arbitrariness of $R(y)(a)$, it follows from (5.39)
that $\varepsilon_r\to 0$ as $r\to \infty$.
In addition, since $\phi_1(\cdot,0),\ldots,\phi_{2n}(\cdot,0)$ satisfy (5.35)-(5.41), it follows from (5.39) and (5.41) that\vspace{-0.1cm}
$$\begin{array}{rrll}
&&\alpha_r=\max\limits_{1\leq i\leq 2n}
\{R(\phi_i)^*(a,0)D_rR(\phi_i)(a,0)\}\\
&\leq &\varepsilon_r \max\limits_{1\leq i\leq 2n}\{\|R(\phi_i)(a,0)\|_1^2\}\\
&\leq &\varepsilon_r(\|E(a)\|_1^2+\|E(a)B(a)\|_1^2+n),
\end{array}                                                                             \vspace{-0.1cm}\eqno (5.43)$$
in which $\Phi(a,0)=(\phi_1,\cdots,\phi_{2n})(a,0)=I_{2n}$ have been used.
Inserting it and $m_r=n_r=0$ into (5.34), we get that (5.33) holds.
The proof is complete. \medskip

\noindent\textit{{\bf Theorem 5.5.} Assume that $\mathscr{L}$ is in l.c.c. at $t=+\infty$.
Then, for each $k\in\Upsilon$, there exists an $r'_k\geq r_k$,
where $r_k$ is specified in Theorem 5.4, such that for all $r\geq r'_k$,\vspace{-0.1cm}
$$|\lambda_{k}^{(r)}-\lambda_{k}|\leq \frac{|\lambda_{k}^{(r)}|^2e_r}{1-|\lambda_{k}^{(r)}|e_r}, \vspace{-0.1cm}\eqno (5.44)$$
$$|\lambda_{k}^{(r)}-\lambda_{k}|\leq \frac{|\lambda_{k}|^2e_r}{1-|\lambda_{k}|e_r},           \vspace{-0.1cm}\eqno (5.45)$$
where $e_r$ denotes the number on the right-hand side in $(5.33)$.}
\medskip

\noindent{\bf Proof.}
For each $k\in\Upsilon$, $\lambda_{k}$ and $\lambda_{k}^{(r)}$ have the same sign for sufficiently large $r$.
In view of $\lambda_k=-1/{\mu_{k}}$ and $\lambda_{k}^{(r)}=-1/{\mu_{k}^{(r)}}$,
it follows from (5.33) that for each $k\in\Upsilon$,\vspace{-0.1cm}
$$|\frac{1}{\lambda_{k}^{(r)}}-\frac{1}{\lambda_{k}}|\leq e_r,\;\;r\geq r_k,\vspace{-0.1cm}$$
which yields that\vspace{-0.1cm}
$$|\lambda_{k}^{(r)}-\lambda_{k}|\leq e_r|\lambda_{k}^{(r)}||\lambda_{k}|,\;\;r\geq r_k. \vspace{-0.1cm} \eqno (5.46)$$
Thus,\vspace{-0.1cm}
$$|\lambda_{k}|=|\lambda_{k}+\lambda_{k}^{(r)}-\lambda_{k}^{(r)}|\leq
|\lambda_{k}-\lambda_{k}^{(r)}|+|\lambda_{k}^{(r)}|\leq
e_r|\lambda_{k}^{(r)}||\lambda_{k}|+|\lambda_{k}^{(r)}|,\;r\geq r_k,\vspace{-0.1cm}$$
which implies that\vspace{-0.1cm}
$$|\lambda_{k}|(1-|\lambda_{k}^{(r)}|e_r)\leq |\lambda_{k}^{(r)}|.                       \vspace{-0.1cm}\eqno (5.47)$$
By Theorem 5.4 and Proposition 5.4, there exists an $r'_k\geq r_k$ such that
$1-|\lambda_{k}^{(r)}|e_r>0$. Hence, it follows from (5.46) and (5.47)
that (5.44) holds. With a similar argument,
one can show that (5.45) holds.
This completes the proof.\medskip

\section{Spectral approximation in the limit point and intermediate cases}
\medskip

Now, we study the regular approximation of spectra of (1.1) in the case that
$\mathscr{L}$ is either in l.p.c. or the intermediate case at $t=+\infty$, namely, $n\leq d<2n$.
In each case, we show that
$\{H_{1,r}\}_{r=1}^\infty$ is spectrally inclusive for any given self-adjoint subspace extension $H_1$.
We always assume that $(\mathbf{A}_1)-(\mathbf{A}_3)$ hold when $\mathscr{L}$ is in l.p.c. at $t=+\infty$
and $(\mathbf{A}_1)-(\mathbf{A}_4)$ hold when $\mathscr{L}$ is in the intermediate case at $t=+\infty$.
\medskip

\noindent\textit{{\bf Theorem 6.1.} Assume that $n\leq d<2n$. Let $H_1$ be any fixed SSE of $H_0$, and $H_{1,r}$ the
induced regular SSE of $H_1$ on $\mathcal{I}_r$, where
$H_1$ and $H_{1,r}$ are determined by (3.9) and (3.12), respectively, when $d=n$ (l.p.c.),
and they are determined by (3.16) and (3.17), respectively, when $n<d<2n$ (the intermediate case).
And let $H_{1,r}'$ be defined by (4.9). Then
\begin{itemize}\vspace{-0.2cm}
\item[{\rm (i)}]  $\{H_{1,r}'\}$ is SRC to $H_1$;\vspace{-0.2cm}
\item[{\rm (ii)}] $\{H_{1,r}\}$ is spectrally inclusive for $H_1$ if $0\not\in \si(H_1)$.
\end{itemize}}

\noindent{\bf Proof.} The main idea of the proof is similar to that of Theorem 5.1, where the core $C(H_1)$
of $H_1$ in (5.3) is replaced by\vspace{-0.2cm}
$$C(H_1)=H_{00}\dotplus L\{\beta_1,\ldots,\beta_d\}, \vspace{-0.2cm}$$
where $\{\beta_i=(\omega_i,\tilde{\tau}_i)\}^d_{i=1}$ is a GKN-set for $\{H_0,H_0^*\}$ and $\omega_i, 1\leq i\leq d$, is defined by (3.8)
and (3.15) when $d=n$ (l.p.c.) and $n<d<2n$ (the intermediate case), separately.
So its details are omitted. The proof is complete.\medskip

\noindent{\bf Remark 6.1.}  In the case that $\mathscr{L}$ is in l.p.c. at $t=+\infty$,
the sequence of induced regular self-adjoint subspace extensions $\{H_{1,r}\}$ is
spectrally inclusive for $H_1$, but not spectrally exact for $H_1$ in general.
For a counterexample, the reader is referred to [19, Example 3.1].
\medskip

\section*{ Acknowledgements }

\setcounter{equation}{0}

The authors would like to thank Mr. H. Zhu for his valuable suggestions for Section 4.

\end{document}